\documentclass{article}

\usepackage{amscd}      
\usepackage{amssymb}     
\usepackage{amsmath}     
\usepackage{xypic}      
\LaTeXdiagrams          
\usepackage[all,v2]{xy}      
\xyoption{2cell}
\UseAllTwocells
\xyoption{frame}
\CompileMatrices
\allowdisplaybreaks[4] 

\usepackage{theorem}

\newtheorem{prop}{Proposition}[section]  
\newtheorem{lem}[prop]{Lemma}

\newtheorem{cor}[prop]{Corollary}
\newtheorem{them}[prop]{Theorem}

\theorembodyfont{\upshape}

\newtheorem{defn}[prop]{Definition}

\newtheorem{numrmk}[prop]{Remark}
\newtheorem{numnote}[prop]{Note}
\newtheorem{numex}[prop]{Example}

\newtheorem{con}[prop]{Construction}

\newtheorem{note}{Note}

\newtheorem{example}{Example}

\newtheorem{rmk}{Remark}

\newenvironment{pf}{\begin{trivlist}\item[]{\sc Proof.}}%
            {\nolinebreak $\Box$ \end{trivlist}}

\newcommand{\noprint}[1]{}

\newcommand{\Section}[1]{\newpage\section{#1}}

\renewcommand{\tilde}{\widetilde}

\newcommand{\qed}{{\nolinebreak $\,\,\Box$}}

\newcommand{\upst}{^{\ast}}
\newcommand{\op}{{{\rm op}}}
\newcommand{\qf}{{\rm qf}}
\newcommand{\f}{{\rm f}}
\newcommand{\perf}{{\rm pf}}
\newcommand{\upsh}{^{!}}
\newcommand{\lst}{_{\ast}}

\newcommand{\lcom}{_{\scriptscriptstyle\bullet}}

\newcommand{\argument}{{{\,\cdot\,}}}

\newcommand{\XX}{{\mathfrak X}}
\renewcommand{\AA}{{\mathfrak A}}

\renewcommand{\SS}{{\mathfrak S}}
\newcommand{\TT}{{\mathfrak T}}

\newcommand{\YY}{{\mathfrak Y}}      
\newcommand{\UU}{{\mathfrak U}}      
      
\newcommand{\ZZ}{{\mathfrak Z}}

\newcommand{\FF}{{\mathfrak F}}         
\newcommand{\GG}{{\mathfrak G}}

\newcommand{\RR}{{\mathfrak R}}
\newcommand{\HH}{{\mathfrak H}}
\newcommand{\VV}{{\mathfrak V}}

\newcommand{\aaa}{{\mathbb A}}

\renewcommand{\O}{{\cal O}}

\newcommand{\cC}{{\cal C}}

\newcommand{\uU}{{\cal U}}

\newcommand{\del}{\partial}

\newcommand{\resto}{{\,|\,}}
\newcommand{\st}{\mathrel{\mid}}

\newcommand{\rk}{\mathop{\rm rk}}

\newcommand{\im}{\mathop{\rm im}}
\newcommand{\ob}{\mathop{\rm ob}}

\newcommand{\cok}{\mathop{\rm cok}}
\newcommand{\amp}{\mathop{\rm amp}}

\newcommand{\spec}{\mathop{\rm Spec}\nolimits}
\newcommand{\dgspec}{\mathop{\mathbb S \rm pec}\nolimits}

\newcommand{\id}{\mathop{\rm id}\nolimits}
\newcommand{\Hom}{\mathop{\rm Hom}\nolimits}

\newcommand{\bhom}{\mathop{\rm\bf Hom}\nolimits}

\newcommand{\shom}{\mathop{\rm Hom}\nolimits^{\scriptscriptstyle\Delta}}

\newcommand{\Homu}{\mathop{\underline{\rm Hom}}\nolimits}
\newcommand{\Deru}{\mathop{\underline{\rm Der}}\nolimits}

\newcommand{\Aut}{\mathop{\rm Aut}\nolimits}

\newcommand{\AAut}{\mathop{\mathfrak Aut}\nolimits}

\newcommand{\injectlim}{\mathop{\lim\limits_{\textstyle\longrightarrow}}\limits}
\newcommand{\projectlim}{\mathop{{\lim\limits_{\textstyle\longleftarrow}}}\limits}
\newcommand{\comp}{\mathbin{{\scriptstyle\circ}}}
\newcommand{\ol}{\overline}
\newcommand{\ul}{\underline}

\newcommand{\longiso}{\stackrel{\textstyle\sim}{\longrightarrow}}

\newcommand{\doublearrowstack}[2]%
                      {{{{\scriptstyle#1}\atop{\textstyle\longrightarrow}}\atop{{\textstyle\longrightarrow}\atop{\scriptstyle#2}}}}
\newcommand{\rightleftarrowstack}[2]%
                      {{{{\scriptstyle#1}\atop{\textstyle\longrightarrow}}\atop{{\textstyle\longleftarrow}\atop{\scriptstyle#2}}}}
\newcommand{\leftrightarrowstack}[2]%
                      {{{{\scriptstyle#1}\atop{\textstyle\longleftarrow}}\atop{{\textstyle\longrightarrow}\atop{\scriptstyle#2}}}}

\newcommand{\ses}[5]%
{0\longrightarrow#1\stackrel{#2}{ \longrightarrow}#3\stackrel{#4}{
\longrightarrow}#5\longrightarrow0}

\newcommand{\dt}[6]%
{#1\stackrel{#2}{longrightarrow}#3 \stackrel{#4}{\longrightarrow}#5
\stackrel{#6}{\longrightarrow} #1[1]}  
 
\newcommand{\cat}[1]%
{(\mbox{\rm #1})}

\title{{\bf \LARGE Differential Graded Schemes~II:}  \\
{\bf \Large  The 2-category of differential graded schemes}}
\author{K. Behrend} 
\date{December 16, 2002}

\begin{document}

\sloppy
\maketitle

\begin{abstract}
We construct a 2-category of differential graded schemes.  The local
affine models in this theory are differential graded algebras, which
are graded commutative with unit over a field of characteristic zero,
are concentrated in non-positive degrees and have perfect cotangent
complex.  Quasi-isomorphic differential graded algebras give rise to
2-isomorphic differential graded schemes and a differential graded
algebra can be recovered up to quasi-isomorphism from the differential
graded scheme it defines. Differential graded schemes can be glued
with respect to an \'etale topology and fibered products of
differential graded schemes correspond on the algebra level to derived
tensor products.
\end{abstract}

\tableofcontents

\newpage
\section*{Introduction}
\addcontentsline{toc}{section}{Introduction}

The goal of this paper is to define a useful notion of differential
graded scheme. This is done with the
following criteria in mind:

(i)  Differential graded schemes can be glued from local data.
Quasi-isomorphisms are considered to be isomorphisms for the purposes
of gluing.

(ii) Every differential graded scheme locally determines a differential
graded algebra up to quasi-isomorphism, the local affine
coordinate ring.

(iii) Fibered products of differential graded schemes exist and are
given locally on the level of affine coordinate rings as derived
tensor products.

(iv) Differential graded schemes form a category.  Thus moduli spaces
such as the ones constructed in \cite{DerQuot} solve  universal
mapping problems in this category and so their differential graded
structure is determined entirely by a universal mapping property.

It turns out that these requirements cannot be met by a usual
category.  Some higher categorical structure is needed.  Our central
message is that the simplest of all higher categorical structures,
namely that of {\em 2-category}, is sufficient for  a satisfying
theory.

This is somewhat of a surprise, because differential graded algebras
form a {\em simplicial category}, which is a version of (weak)
infinity category.  Passing to 2-categories is achieved by a process
of {\em truncation}, which leads, by its nature, to loss of
information.  The fact that the lost information was not necessary for
the purposes of geometry is rather subtle.  It is the content of the
results we present under the heading of {\em Descent Theory}.

\subsubsection{Overview of the construction}

We start with a suitable 2-category of differential graded algebras.
This is the 2-category of {\em perfect resolving algebras}, which we
denote by $\RR_\perf$. A perfect resolving algebra is a differential
graded algebra concentrated in non-positive degrees (the differential
has degree $+1$), such that the underlying graded algebra is free
(commutative with 1, over a field $k$ of characteristic $0$) on
finitely many generators in each degree, and such that its complex of
differentials is {\em perfect}.  For a detailed study of perfect
resolving algebras, see \cite{dgsI}. 

The perfect resolving algebras form a full subcategory of the
differential graded algebras, which form a simplicial closed model
category. Thus, given any two perfect resolving algebras $B$, $A$, the
set of morphisms from $B$ to $A$ is a simplicial set
$\shom(B,A)$. Because perfect resolving algebras are both fibrant and
cofibrant, the simplicial set $\shom(B,A)$ is fibrant, i.e., it has
the Kan extension property.  Thus the fundamental groupoid of
$\shom(B,A)$ exists and we define
$$\bhom(B,A)=\Pi_1\shom(B,A)\,.$$
With this definition of hom-groupoid, the perfect resolving algebras
form the 2-category $\RR_\perf$. Let us note that two perfect
resolving algebras are isomorphic in $\RR_\perf$ (in the `relaxed',
2-categorical sense), if and only if they are quasi-isomorphic.  This
is because quasi-isomorphisms are the {\em weak equivalences }in the
closed model category of differential graded algebras.

The 2-category $\RR_\perf$ serves as the category of affine
coordinate rings of affine differential graded schemes. To construct
differential graded schemes over $\RR_\perf$, we
imitate the usual construction of algebraic spaces over the category
of $k$-algebras (of finite type, to keep the analogy with
$\RR_\perf$). 

Thus, the first step is to pass to the opposite category of
$\RR_\perf$, which we denote $\SS$.  Then we introduce a Grothendieck
topology on the 2-category $\SS$.  The usual \'etale topology on the
category of affine $k$-schemes of finite type has an analogue in
$\SS$, called, not surprisingly, the {\em \'etale }topology on $\SS$.
The necessary facts about \'etale morphisms between perfect resolving
algebras are proved in~\cite{dgsI}.

As soon as we have a 2-category with a Grothendieck topology, we have
the category of {\em sheaves }over it.  We call a sheaf over $\SS$ a
{\em differential graded sheaf}. {\em Differential graded schemes }are
then defined to be differential graded sheaves satisfying an extra
condition (see below).

Since we are working over a
2-category, the notion of sheaf resembles more closely the usual
concept of {\em stack}, rather than the usual concept of sheaf.  In
fact, a sheaf over $\SS$ is defined to be a {\em 2-category fibered in
groupoids }over $\SS$.  It is required to satisfy sheaf axioms, which
are direct adaptations of the usual stack axioms.  Thus, on a certain
formal level, our theory of differential graded schemes resembles the
usual theory of algebraic stacks.  

There is one essential difference: there is no 1-category which
generates the 2-category of differential graded schemes in the same
way that the 1-category of usual schemes generates the 2-category of
algebraic stacks. The local affine models for differential graded
schemes already form a 2-category, in contrast to the local affine
models for algebraic stacks, which form only a 1-category.

A key ingredient in the construction of usual algebraic spaces is {\em
descent theory}. By this we mean two results: descent for morphisms
and descent for algebras. Descent for morphisms says that the
contravariant functor 
\begin{equation}\label{specr}
\spec (R):(\text{finite type affine $k$-schemes})\longrightarrow
(\text{sets})
\end{equation}
represented by the finite type $k$-algebra $R$, is a sheaf. Thus we
obtain a contravariant functor 
$$\spec:(\text{finite type $k$-algebras})\longrightarrow
\big(\text{sheaves on $(\text{finite type affine
$k$-schemes})$}\big)\,.$$ 
By Yoneda's lemma, it is fully faithful.  Without descent theory, we
would have to pass from $\spec R$ to the associated sheaf, which would
destroy the fully faithful property of $\spec$.  As a consequence, we
would not be able to reconstruct the finite type $k$-algebra $R$ from
the sheaf (and hence the algebraic space) associated to it.

Thus, in view of our requirement (ii) on differential graded schemes,
a result on descent for morphisms in $\SS$ is essential.  It says that
the 2-category fibered in groupoids over $\SS$
$$\dgspec(B)$$
represented by the perfect resolving algebra $B$ is a sheaf.
From  the `lax functor' point of view, this 2-category fibered in
groupoids $\dgspec( B)$ may be considered as a contravariant 2-functor
$$\dgspec( B):\SS\longrightarrow(\text{groupoids})\,,$$
making the analogy with (\ref{specr}) more apparent. 

We obtain a contravariant 2-functor
$$\dgspec: \RR_\perf\longrightarrow (\text{sheaves on $\SS$})\,.$$ 
By Yoneda's lemma for 2-categories, it is fully faithful.  Again, the
key point is that there is no need to pass to an associated sheaf, and
so there is no information loss when passing from a perfect resolving
algebra to the sheaf on $\SS$ it gives rise to. This means that our
requirement (ii), above, is fulfilled.

Usual descent for algebras can be formulated as follows.  Call a
morphism of sheaves $f:X\to Y$ on the category of finite type
$k$-algebras {\em affine}, if for every morphism $\spec R\to Y$, the
fibered product $X\times_Y\spec R$ is isomorphic to $\spec S$, for
some finite type $k$-algebra $S$. Descent for algebras says that for
$f:X\to Y$ to be affine it suffices to have an \'etale cover $\spec
R_i\to Y$ of $Y$, such that for every $i$, the fibered product
$X\times_Y\spec R_i$ is isomorphic to $\spec S_i$, for some finite
type $k$-algebra $S_i$.  We abbreviate this property by saying that
`affine' is a {\em local }property for morphisms between sheaves.  When
developing the theory of algebraic spaces from the theory of affine
schemes, this fact is essential.  Without it, it would be impossible
to ever check that any given morphism is affine.  By extension, it
would be impossible to ever prove that a given sheaf is an algebraic
space.

Thus we prove an analogue of descent for algebras in the 2-category
$\SS$.  Once this is done, we define differential graded schemes in
three steps:

$\bullet$ An {\em affine }differential graded scheme is a differential
graded sheaf, 2-isomorphic to  $\dgspec B$, for some perfect resolving
algebra $B$.

$\bullet$ An {\em affine morphism }of differential graded sheaves is
a morphism, whose base change to an affine differential graded scheme
always gives rise to an affine differential graded scheme.  For affine
morphisms, the property of being {\em \'etale }makes sense.

$\bullet$ A {\em differential graded scheme }is a differential graded
sheaf $\XX$, which can be covered by affine \'etale morphisms $\dgspec
B_i\to\XX$.  Thus, a differential graded scheme is \'etale locally
affine. 

By the local nature of this definition, it is clear that it satisfies
our criterion~(i), above.  A more detailed study of the gluing
properties of differential graded schemes is the content
of~\cite{dgsIII}.  There we will prove, for example, that every local
complete intersection scheme can be considered as a differential
graded scheme.  Requirement~(iv) is also clearly fulfilled: morphisms
between differential graded schemes are just morphism of differential
graded sheaves.  Because of the truncation procedure involved in our
construction, Property~(iii) is somewhat non-trivial.  It will be
dealt with in~\cite{dgsIV}.  For the present purposes is sufficient to
have base changes by \'etale morphisms (see
Proposition~\ref{base.et.r}). 

Recall (Theorem~\ref{loc.fin} in~\cite{dgsI}), that every perfect
resolving algebra is locally finite. Thus every differential graded
scheme can be glued using only finite resolving algebras. In other
words, the 2-category of affine differential graded schemes associated
to finite resolving algebras generates the 2-category of differential
graded schemes.  Thus we could base our theory entirely on finite
resolving algebras instead of perfect resolving algebras. On the other
hand, the descent result for algebras fails in the context of finite
resolving algebras. Hence finite resolving algebras give rise
differential graded schemes which are somewhat too local, to be
considered as the class of {\em all }affine differential graded
schemes.  

Because of these observations, one might speculate that it
should be possible to develop the theory of differential graded
schemes without descent for algebras.  On the other hand, the
development is greatly simplified by its use.

\subsubsection{Outline of the paper}

In Section~\ref{sec.etatop}, we start by reviewing a few basic facts
about 2-categories. Then we define presheaves over 
2-categories. We introduce the  notion of 
Grothendieck topology on a 2-category and define sheaves on a
2-category endowed with a Grothendieck topology.  

We proceed to introduce the 2-category of resolving algebras $\RR$,
together with its subcategories of {\em quasi-finite}, {\em perfect
}and {\em finite }resolving algebras, $\RR_\qf$, $\RR_\perf$ and
$\RR_\f$.  We obtain a base 2-category $\SS$ by passing to the
opposite 2-category of any of $\RR_\qf$, $\RR_\perf$ or $\RR_\f$. We
prove that base changes by \'etale morphisms exist in $\SS$.

Finally, we introduce the \'etale topology on $\SS$ and define
differential graded sheaves as sheaves on $\SS$.

Section~\ref{sec.desc} contains our results on descent theory.  There
is, first of all, a theorem on descent for morphisms:
Theorem~\ref{Descent} and its Corollary~\ref{sheaf.prop}.  It holds
for both finite and perfect resolving algebras.  This result is
really the technical heart of the whole theory, because it justifies
using 2-categories for differential graded schemes.  Without it, one
would have to consider some type of infinity category (as is done
in~\cite{ToenVezz}). To prove our descent theorem, we use the main
result of~\cite{dgsI}, on `linearization of homotopy groups'.  It says
that for every $\ell>0$ there exists a canonical bijection
$$\Xi_\ell:h^{-\ell}\Deru(B,A)\longrightarrow\pi_\ell\shom(B,A)\,,$$
where $\Deru(B,A)$ is the differential graded $A$-module of
(internal) derivations $D:B\to A$.  Thus $\pi_\ell\shom(B,A)$ can be
given the structure of $h^0(A)$-module, which suffices to reduce
descent for morphisms in $\SS$ to usual descent for morphisms in the
category of finite type $k$-algebras. 

To formulate our theorem on descent for algebras, we need to introduce
a special class of {\em gluing data }in $\SS$.  We leave the general
theory of gluing data in $\SS$ to~\cite{dgsIII}. Here we only require
{\em relative }gluing data, which have a much simpler structure. These
gluing data can be conveniently pictured as diagrams in the shape of
truncated hypercubes in $\SS$.  We prove that relative gluing data can
be {\em  strictified}, which means that 2-arrows appearing in the
squares of the hypercube can be replaced by identity 2-arrows. 

Once these preliminaries are dispensed with, we proceed to prove our
theorem on descent for algebras, Theorem~\ref{desalg}. As mentioned
above, it is the basis for the theory of affine morphisms of
differential graded sheaves and is therefore also an integral part of
the definition of differential graded scheme.  

Section~\ref{sec.schemes} contains this definition, as outlined
above. At this point, we refrain from going much further than the bare
definition of differential graded scheme. 

Instead, in Section~\ref{sec.oneinv}, we proceed to construct the
basic `1-categorical invariants' of a differential graded scheme
$\XX$.  All of these are  sheaves (or complexes of sheaves) on the
1-category associated to the 2-category underlying $\XX$.  An object
of this 1-category may be thought of as an isomorphism class of
morphisms 
$$\dgspec A\to \XX\,.$$
There are
first of all the {\em higher structure sheaves}. These associate to
$\dgspec A\to \XX$ the $h^0(A)$-module $h^i(A)$, for $i\leq0$. 

Then there are the {\em higher tangent sheaves
}$h^{\ell}(\Theta_\XX)$.  If $\XX=\dgspec B$
is affine, they associate to $\dgspec A\to \dgspec B$ the
$h^0(A)$-module $h^\ell\Deru(B,A)$, for various $\ell$. 

Next there are the {\em homotopy sheaves
}$\pi_\ell(\XX)$, for $\ell>0$. Again,
let us just say here that in the affine case $\XX=\dgspec B$, they are
given by associating to $\dgspec A\to \dgspec B$ the group
$\pi_\ell\shom(B,A)$.  These are, in fact, sheaves, by
our results on descent theory. Moreover, they coincide with certain of
the higher tangent sheaves: there exists a canonical isomorphism of
sheaves
$$\Xi_\ell:h^{-\ell}(\Theta_\XX)\longiso\pi_\ell(\XX)\,,$$
for all $\ell>0$. 

To put the homotopy sheaves $\pi_\ell(\XX)$ into context, let us make
a few general remarks.  Let $\cC$ be a simplicial closed model category with
homotopy category $Ho(\cC)$, say in its incarnation as category of
fibrant-cofibrant objects with simplicial homotopy classes of maps as
morphisms.  For every $\ell\geq0$ and every morphism $f:X\to Y$
in $Ho(\cC)$ we let $\pi_\ell(X/Y)$ be the presheaf on 
$Ho(\cC)/X$ defined by $\pi_\ell(X/Y)(U)=\pi_\ell\shom_Y(U,X)$, for all
$U\to X$ in $Ho(\cC)/X$.  Here $\shom_Y(U,X)$ denotes the fiber of
$\shom(U,X)\to\shom(U,Y)$.  For $Y=\ast$, we obtain the presheaf
$\pi_\ell(X)$.  There is a long exact sequence of presheaves of
pointed sets on $Ho(X)$
\begin{multline*}
\ldots\longrightarrow \pi_\ell(X/Y) \longrightarrow \pi_\ell(X)
\longrightarrow f^{-1}\pi_\ell(Y)\stackrel{\del}{\longrightarrow}
\pi_{\ell-1}(X/Y)\longrightarrow\ldots\\
\ldots\longrightarrow f^{-1}\pi_1(Y)\stackrel{\del}{\longrightarrow}
\pi_0(X/Y) \longrightarrow \pi_0(X)\longrightarrow \pi_0(Y)\,.
\end{multline*}

Formally, our homotopy sheaves $\pi_\ell(\XX)$, of which there also
exist relative versions $\pi_\ell(\XX/\YY)$, behave somewhat  as if
they were obtained, as above, from a simplicial closed model category
structure underlying the 2-category of differential graded schemes.
In particular, they also fit into a long exact sequence.  We do not
know if there exists such a simplicial closed model category structure
underlying the 2-category of simplicial schemes, but we find it quite
likely.

There is a special feature in our case of differential graded
schemes: the analogue of the homotopy category $Ho(\cC)/X$ is the
1-category associated to the differential graded scheme $\XX$.  Thus
in our case this homotopy category is endowed with a Grothendieck
topology, with respect to which all the $\pi_\ell(X)$, for $\ell>0$, are
sheaves.  This property does not seem to have a meaningful analogue,
for example, in the simplicial closed model categories of topological
spaces or simplicial sets.  

Finally,we define the cotangent complex of a morphism of differential
graded schemes and the algebraic space associated to a differential
graded scheme. The associated algebraic space is obtained by gluing
the truncations $\spec h^0(B_i)$ of the affine differential graded
schemes $\dgspec B_i$, which cover a differential graded scheme $\XX$.
The cotangent complex of a differential graded scheme gives rise to an
{\em obstruction theory }in the sense of \cite{BF} on the associated
algebraic space.  Thus in the case of perfect amplitude 1, it defines
a {\em virtual fundamental class }on the associated algebraic space.

\subsubsection{Notation}

References of the form I.0.0, refer to Result~0.0 of~\cite{dgsI}.

\subsubsection{Acknowledgements}

I wish to thank the Research Institute for the Mathematical Sciences
in Kyoto, for providing very comfortable working conditions during my
visit in 1999/2000. I would especially like to thank my host,
Professor K. Saito, for the warm hospitality.  Most of the research
presented in this paper was carried out during my sabbatical at RIMS
in 1999/2000. In particular, the discovery of the 
\'etale topology and the discovery of the  utility of 2-categories
date back to that period.

I would also like to thank the organizers of the workshop on Algebraic
Geometry and Integrable Systems related to String Theory at RIMS in
June 2000, especially Professor M.-H. Saito, for giving me the first
opportunity to present this theory of differential graded schemes to
the public.

In addition, I would  like to thank the Mathematical Sciences Research
Institute in Berkeley.  I worked out some of the details of descent
theory during a visit to MSRI in early 2002.

This work was also continuously supported by a Research Grant from the
National Science and Engineering Council of Canada. 

Finally, I would like to thank I. Ciocan-Fontanine, E. Getzler and
M. Kapranov for helpful discussions.

\Section{The \'etale topology}\label{sec.etatop}

\subsection{2-Categories}

All our 2-categories will have invertible 2-morphisms. Thus they are
categories enriched over groupoids.  To fix notation, let us recall
the definition. (See also~\cite{CftWM}). 

\begin{defn}\label{twocat}
A {\bf 2-category }$\SS$ consists of 

(i) a set of objects $\ob\SS$,

(ii) for every pair $U$, $V$ of objects of $\SS$ a groupoid
$\bhom(U,V)$,

(iii) for every triple $U$, $V$, $W$ of objects of $\SS$ a functor
\begin{align}\label{horcomp}
\comp:\bhom(V,W)\times\bhom(U,V)& \longrightarrow \bhom(U,W)\\*
(f,g)&\longmapsto f\comp g\nonumber
\end{align}

(iv) for every object $U$ of $\SS$ an object $\id_U$ of
$\bhom(U,U)$, or rather a functor
\begin{equation*}
\ast \stackrel{\id_U}{\longrightarrow} \bhom(U,U),
\end{equation*}

such that 

(i) the composition $\comp$ is associative, i.e., for four objects
$U$, $V$, $W$, $Z$ of $\SS$ the diagram of functors
$$\begin{diagram}
{\bhom(W,Z)\times\bhom(W,V)\times\bhom(V,U)}\dto\rto &
{\bhom(V,Z)\times\bhom(U,V)}\dto\\
{\bhom(W,Z)\times\bhom(U,W)}\rto & {\bhom(U,Z)}
\end{diagram}$$
commutes (strictly),

(ii) $\id_U$ acts as identity, i.e., the induced diagrams 
$$\xymatrix@C=1pc{
{\bhom(U,V)}\dto\drto & {\bhom(U,V)}\drto\rto &
{\bhom(V,V)\times\bhom(U,V)}\dto\\
{\bhom(U,V)\times\bhom(U,U)}\rto& {\bhom(U,V)} & {\bhom(U,V)}}$$
commute (strictly), for any two objects $U$, $V$.

Composition in $\bhom(U,V)$ is called {\bf vertical }composition, the
operation $\comp$ of~(\ref{horcomp}) is called {\bf horizontal
}composition. Vertical composition we shall denote by
$\alpha\cdot\beta$.

Objects of $\bhom(U,V)$ are called {\bf 1-morphisms }(of $\SS$),
morphisms in $\bhom(U,V)$ are called {\bf 2-morphisms }(of $\SS$). We
also use the words {\em 2-isomorphism }and {\em 2-arrow }instead of
2-morphism.  An identity 2-arrow is also called a {\em strictly
commutative diagram}. 
\end{defn}

Every set $X$ is a category by taking the elements of $X$ as objects
and admitting only identity morphisms.  Every category $\cC$ is a 2-category,
by considering the Hom-set $\Hom(A,B)$ as a category, for any two
objects $A$, $B$ of $\cC$. 

Given a 2-category $\SS$, the objects of $\SS$ together with the
1-morphisms and horizontal composition form a 1-category, the {\bf
underlying }1-category of $\SS$.  Replacing $\bhom(A,B)$ by its set of
isomorphism classes, we obtain another 1-category, the 1-category {\bf
associated }to $\SS$, which we denote by $\ol\SS$.  There is a
canonical functor from the underlying 1-category of $\SS$ to $\ol\SS$.

Compatibilities between 2-morphisms can often be phrased conveniently by
saying that certain `2-spheres' commute. This means
that the objects involved should be considered as vertices, 1-morphisms as
edges and 2-morphisms as faces of a `triangulation' of a topological
2-sphere. There should always be one `source object', having only
1-morphisms emanating from it, and one `target object', which has no
1-morphism emanating from it. Then all the different directed paths
from the source to the target object can be considered as vertices of
a commutative polygon of 2-arrows, i.e., a commutative polygon for the
vertical composition.
Often we project such a 2-sphere stereographically
onto the plane, so that we get a flat diagram, whose exterior should be
considered as a 2-cell, even if it is not labelled as such.

If the distinction is important, we say that such a 2-sphere
{\em 2-commutes}.  This is contrast to a 2-sphere all of whose faces
are identity 2-arrows (in other words, are commutative diagrams in the
underlying 1-category), which we call {\em strictly }commutative.  If
all faces are strictly commutative, a 2-sphere is automatically
2-commutative.

\subsubsection{Isomorphisms and fibered products in 2-categories}

\begin{defn}
For lack of better terminology, we will call a morphism of groupoids
$f:X\to Y$ {\bf categorically \'etale}, if for every object $x$ of $X$
the induced group homomorphism
$\Aut_X(x)\to\Aut_Y\big(f(x)\big)$ is bijective. 
\end{defn}

Let $\SS$ be a 2-category.

\begin{defn}
We call a 1-morphism $A\to B$ in $\SS$ {\bf faithful} ({\bf
categorically \'etale}, a {\bf monomorphism}), if for
every object $S$ the 
induced morphism of groupoids $\bhom(S,A)\to\bhom(S,B)$ is faithful
(categorically \'etale, fully faithful). 
\end{defn}

Thus we have the implications:
$$\text{\rm monomorphism}\Longrightarrow \text{categorically
\'etale}\Longrightarrow \text{faithful}\,.$$

\begin{defn}
Let $f:A\to B$ be a 1-morphism in a 2-category $\SS$.  A {\bf
2-inverse }of $f$ is given by the data $(g,\phi,\psi)$, where $g:B\to
A$ is a 1-morphism and $\phi:\id_A\Rightarrow g\comp f$ and
$\psi: f\comp g\Rightarrow\id_B$ are 2-arrows, such that the two
diagrams
\vskip-.3cm
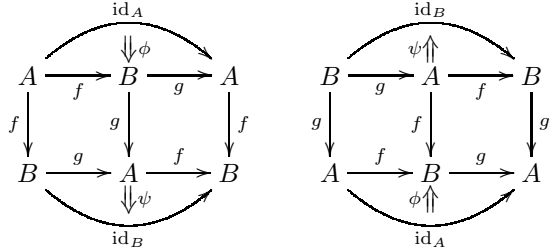
$$\begin{diagram}
{A\dto_f\rto_f\rruppertwocell^{\id_A}<8>{_\phi}} & 
B\dto_g\rto_g & A\dto^f &
{B\dto_g\rto_g\rruppertwocell^{\id_B}<8>{^\psi}} & 
A\dto_f\rto_f & B\dto^g\\
{B\rto^g\rrlowertwocell_{\id_B}<-8>{_\psi}} & A\rto^f & B &
{A\rto^f\rrlowertwocell_{\id_A}<-8>{^\phi}} & B\rto^g & A
\end{diagram}$$
\vskip-.5cm
\noindent commute.
\end{defn}

The inverse of $f:A\to B$ is unique up to a unique 2-isomorphism in
the following sense.  If $(g',\phi',\psi')$ is another inverse to $f$,
then there exists a  unique 2-isomorphism $\theta:g\to g'$, such that
the two diagrams
\vskip-.8cm
$$\begin{diagram}
A\rto^f\rruppertwocell<12>^{\id_A}{_\phi}
\rrlowertwocell<-12>_{\id_A}{^\phi'}&
B\rtwocell^g_{g'}{_\theta}& A &
B\rtwocell^g_{g'}{_\theta}\rruppertwocell<12>^{\id_B}{^\psi}
\rrlowertwocell<-12>_{\id_B}{_\psi'}&
A\rto^f& B
\end{diagram}$$
\vskip-.8cm
\noindent commute.

\begin{defn}
We call a 1-morphism $f:A\to B$ in a 2-category {\bf 2-invertible},
an {\bf equivalence} or even an {\bf isomorphism}, if it admits a
2-inverse.  We hope that there will be no confusion with the  term
{\em 2-isomorphism}, which stands for the  2-arrows in $\SS$.  

By the same token, we call two objects $A$ and $B$ {\bf isomorphic},
if there exists a 2-invertible morphism $f:A\to B$.
\end{defn}

The following is a very useful criterion by which to recognize
equivalences:

\begin{prop}
Let $f:A\to B$ be a 1-morphism in $\SS$. If there  exists a 1-morphism
$g:B\to A$ such that $\id_A\cong g\comp f$ and $\id_B\cong f\comp g$,
then $f$ is an equivalence.\qed
\end{prop}

\begin{cor}
A 1-morphism $A\to B$ in $\SS$ is an isomorphism, if and only if
$\bhom(S,A)\to \bhom(S,B)$ is an equivalence of groupoids, for all
objects $S$.
\end{cor}

\begin{defn}
A diagram 
\begin{equation}\label{strongweak} 
{\begin{diagram}
{W\dto\rto\drtwocell\omit{^}} &{ B\dto}\\
{A\rto} &{ Z}
\end{diagram}} \end{equation} 
in a 2-category $\SS$ is called a {\bf fibered product}, if for
any object $S$ of $\SS$ the functor
\begin{equation}\label{wkfiber}
\bhom(S,W)\longrightarrow \bhom(S,A)\times_{\bhom(S,Z)}\bhom(S,B)
\end{equation}
is an equivalence of groupoids. Here, the
fibered product in (\ref{wkfiber}) is the usual fibered product of
groupoids. 

We also say that Diagram~(\ref{strongweak}) is {\bf 2-cartesian}. 
\end{defn}

A {\em 1-cartesian }diagram in the 2-category $\SS$ is a strictly
commutative square, which is cartesian in the underlying 1-category. 

\subsubsection{2-Functors and natural 2-transformations}

\begin{defn}\label{twofunctor}
Let $\SS$ and $\TT$ be 2-categories.  A {\bf 2-functor }$f:\SS\to\TT$
consists of 

(i) a map  $f:\ob\SS\to\ob\TT$,

(ii) for every pair $U$, $V$ of objects of $\SS$ a functor 
$$f:\bhom(U,V)\longrightarrow\bhom(f(U),f(V)),$$

\noindent such that 

(iii) for every object $U$ of $\SS$ the diagram of functors
$$\begin{diagram}
\ast\dto_{\id_U}\drto^{\id_{f(U)}} & \\
\bhom(U,U)\rto^-f & \bhom(f(U),f(U))
\end{diagram}$$
commutes (strictly),

(iv) for every triple $U$, $V$, $W$ of objects of $\SS$
$$\begin{diagram}
{\bhom(V,W)\times\bhom(U,V)\dto_{\comp}\rto^-f} &
{\bhom(f(V),f(W))\times\bhom(f(U),f(V))\dto_{\comp}}\\
{\bhom(U,W)\rto^-f} &{\bhom(f(U),f(W))}
\end{diagram}$$
is a (strictly) commutative diagram of functors.
\end{defn}

\begin{defn}
Consider a 2-functor $f:\SS\to\TT$.

(i) $f$ is called {\bf fully faithful}, if for any
two objects $U$, $V$ of $\SS$
the functor $\bhom(U,V)\to\bhom\big(f(U),f(V)\big)$ is an equivalence
of groupoids.

(ii) $f$ is called an {\bf equivalence of 2-categories}, if it is
fully faithful, and for every object $T$ of $\TT$, there exists an
object $S$ of $\SS$, such that $f(S)\cong T$.
\end{defn}

\begin{defn}\label{twotransformation}
Let $f:\SS\to\TT$ and $g:\SS\to\TT$ be two 2-functors between the
2-categories $\SS$ and $\TT$.  A {\bf natural 2-transformation
}$\theta:f\to g$ is given by a functor
$$\ast\stackrel{\theta(U)}{\longrightarrow}\bhom(f(U),g(U)),$$
for every object $U$ of $\SS$, such that 
$$\xymatrix@C-4pc{
{\bhom(U,V)\dto_{g\times\theta(U)}\rto^-{\theta(V)\times f}}&
{\bhom(f(V),g(V))\times \bhom(f(U),f(V))\dto^{\comp}}\\
{\bhom(g(U),g(V))\times\bhom(f(U),g(U))\rto^-{\comp}}&
{\bhom(f(U),g(V))}
}$$
is a commutative diagram of functors, for any two objects $U$, $V$ in
$\SS$.
\end{defn}

\begin{rmk}
2-categories, 2-functors and natural 2-transformations
form a 2-category.
\end{rmk}

\begin{defn}
Let $\SS$ be a 2-category.  We define the {\bf opposite }2-category
$\SS^\op$ as follows:

(i) objects of $\SS^\op$ are the same as the objects of $\SS$,

(ii) given objects $U$ and $V$ of $\SS^\op$, we set
$$\bhom_{\SS^\op}(U,V)=\bhom_{\SS}(V,U)^\op,$$

(iii) horizontal composition in $\SS^\op$ is defined to be the
`opposite' of horizontal composition in $\SS$, up to changing the
order of the arguments,

(iv) identity objects in $\SS^\op$ are the same as in $\SS$.

A {\bf contravariant }2-functor $\SS\to\TT$ is a 2-functor
$\SS^\op\to\TT$.
\end{defn}

Of course, the opposite 2-category is again a 2-category.

\subsection{Presheaves}

We shall now introduce the notion of a {\em presheaf }over a
2-category.  This generalizes the notion of {\em category fibered
in groupoids}, known from the theory of 1-categories (and used in the
theory of stacks).  If the base category is a 2-category, it is 
natural to think of this structure as generalizing the notion of
presheaf over a 1-category.

\begin{defn}\label{defnfibration}
A 2-functor $\pi:\FF\to\SS$ of 2-categories is called a {\bf
presheaf}, if 

(i) for every 1-morphism 
$$\begin{diagram}
V\rto & U\end{diagram}$$
in $\SS$, and every object $x$ of $\FF$ lying over $U$, there exists a
1-morphism 
$$\begin{diagram}
y\rto & x\end{diagram}$$
in $\FF$, lying over $V\to U$,

(ii) for every 2-commutative triangle
\begin{equation}\label{comtri}
{\begin{diagram}
V'\rto^{}\drlowertwocell<0>_{f'}{^<-2>} & V\dto^{f}\\
& U\end{diagram}}
\end{equation}
in $\SS$, and every diagram  of 1-morphisms
$$\begin{diagram}
y'\drto_{e'}& y\dto^{e}\\
& x\end{diagram}$$
in $\FF$, where $e'$ lies over $f'$ and $e$ lies over $f$,
there exists a unique pair $(r,\gamma)$ in $\FF$ such that the triangle
\begin{equation}\label{triup}
{\begin{diagram}
y'\rto^{r}\drlowertwocell<0>_{e'}{^<-2>\gamma} & y\dto^{e}\\
& x\end{diagram}}
\end{equation}
2-commutes in $\FF$ and the triangle~(\ref{triup}) lies over the
triangle~(\ref{comtri}).
\end{defn}

The following special case of Condition~(ii) is worth pointing out:

\begin{lem}
Let $\pi:\FF\to\SS$ be a presheaf.  Then for every 2-morphism
$$\begin{diagram}V\rtwocell^{f'}_{f}{\alpha} & U
\end{diagram}$$
in $\SS$, whose target $f:V\to U$ has been lifted to a 1-morphism 
$$\begin{diagram}y\rlowertwocell<-2>_{e}{\omit} & x
\end{diagram}$$
in $\FF$, there exists a unique 2-morphism 
$$\begin{diagram}y\rtwocell^{e'}_{e}{\gamma} & x
\end{diagram}$$
in $\FF$ with target $e$, which lies over $\alpha$.
\end{lem}
\begin{pf}
Apply the case that $y=x$ and $e$ the identity in Condition~(ii) of
Definition~\ref{defnfibration} to the inverse of $\alpha$. 
\end{pf}

Note that the fibers of a presheaf are groupoids.  We denote the fiber of
$\FF\to\SS$ over the object $U$ of $\SS$ by $\FF_U$.

\begin{defn}
If $\pi:\FF\to\SS$ is a presheaf, and the lifting required in
Condition~(i) of Definition~\ref{defnfibration} is always {\em
unique}, then we call $\FF$ a {\bf presheaf of sets}.  Note that
this makes Condition~(ii) empty.  Moreover, it implies that all fibers of
$\pi$ are sets.
\end{defn}

\begin{defn}
Let $\FF$ and $\GG$ be presheaves over $\SS$.  A {\bf morphism }of
presheaves is a 2-functor
$\FF\to\GG$ commuting with the projections to $\SS$.  A {\bf
2-morphism }of morphisms of presheaves is a 
natural 2-transformations lying over the identity natural
2-transformation. 
\end{defn}

The presheaf morphisms from $\FF$ to $\GG$ form a groupoid
$\bhom_\SS(\FF,\GG)$, the {\em groupoid of presheaf morphisms }from
$\FF$ to $\GG$. 

There is an obvious way to define a composition functor
$$\bhom_\SS(\GG,\HH)\times\bhom_\SS(\FF,\GG)
\longrightarrow\bhom_\SS(\FF,\HH)\,,$$
thus turning the presheaves over $\SS$ into a 2-category. 

Since the fibers of a presheaf over $\SS$ are groupoids, presheaves behave
more like categories than like 2-categories:

\begin{lem}\label{morpre}
Let $\FF$ and $\GG$ be presheaves over $\SS$. A morphism $F:\FF\to\GG$ is a
1-functor of the underlying 1-categories of $\FF$ and $\GG$, commuting with
the underlying 1-functors to the underlying 1-category $\SS$, with the
following property: if 
$$\xymatrix{y\rtwocell^f_g{\theta} & x}$$
is a 2-morphism in $\FF$, lying over
\begin{equation}\label{morpreeq}
\xymatrix{V\rtwocell & U\,,}
\end{equation}
then there exists a 2-morphism
$$\xymatrix{F(y)\rtwocell^{F(f)}_{F(g)}{\eta} & F(x)}$$
lying over (\ref{morpreeq}), also.

Let $F,G:\FF\to\GG$ be 1-morphisms of presheaves over $\SS$. A
2-isomorphism $\theta:F\Rightarrow G$ is a natural 1-transformation between
the underlying 1-functors $F$ and $G$, which maps every object $x$ of $\FF$
to a morphism of $\GG$ lying over $\id_U$ (where $x$ lies over $U$).  The
compatibility with 2-isomorphisms is then automatic. \qed
\end{lem}

\begin{prop}\label{monocrit}
Let $F:\FF\to\GG$ be a morphism of presheaves over $\SS$. The
following are equivalent:

(i)  $F$ is faithful (categorically \'etale, a monomorphism, an
isomorphism)
as a 1-morphism in the 2-category of presheaves over $\SS$,

(ii) for every object $U$ of
$\SS$, the fiber 
functor 
$F_U:\FF_U\to\GG_U$ is faithful (categorically \'etale, fully
faithful, an equivalence).\qed 
\end{prop}

\subsubsection{Relative 2-categories, Yoneda Theory}

Given an object $U$ of a 2-category $\SS$, we define the {\em relative
2-category }$\SS/U$ as follows.

Objects of $\SS/U$ are 1-morphisms $V\to U$ in $\SS$ with target $U$.  Given
two such relative objects $f':V'\to U$ and $f:V\to U$, we define the
groupoid of morphisms $\bhom_{\SS/U}(V',V)$ to have as objects the
2-commutative diagrams
$$\begin{diagram}
V'\rto^{e}\drlowertwocell<0>_{f'}{^<-2>\alpha} & V\dto^{f}\\
& U,\end{diagram}$$
and to have as 2-morphisms from $(e,\alpha)$ to $(g,\beta)$ the
2-morphisms in $\SS$
$$\begin{diagram}V'\rtwocell^{e}_{g}{\theta} & V,
\end{diagram}$$
such that $(f\comp\theta)\cdot\alpha=\beta$, i.e., the diagram of
2-arrows in $\SS$
$$\begin{diagram}
f'\ar@{=>}[r]^{\alpha}\ar@{=>}[dr]_{\beta} & f\comp e
\ar@{=>}[d]^{f\comp \theta}\\ 
&f\comp g\end{diagram}$$
commutes.

Composition in $\bhom_{\SS/U}(V',V)$ is induced from vertical
composition in $\SS$.  We define horizontal composition $\SS/U$
by the formula
$$(e,\alpha)\comp(g,\beta)=\big(e\comp g, (\alpha\comp
g)\cdot\beta\big)\,.$$

By projecting onto the source, more precisely, mapping $V\to U$ to
$V$, $(e,\alpha)$ to $e$ and $\theta$ to $\theta$, we get a 2-functor
$\SS/U\to \SS$.

\begin{lem}
The 2-functor $\SS/U\to\SS$ is a presheaf. \qed
\end{lem}

Abbreviate for an object $U$ of $\SS$ the relative 2-category $\SS/U$ by
$\ul{U}$.

The association $U\to\ul{U}$ defines a 2-functor from $\SS$ to the
2-category of presheaves over $\SS$.  The analogue of Yoneda's lemma in
this context is that this 2-functor is fully faithful:

\begin{prop}[Yoneda's lemma for 2-categories]\label{2Y}
The 2-functor 
\begin{align*}
\SS & \longrightarrow (\text{\rm presheaves/$\SS$})\\*
U & \longmapsto \ul{U}
\end{align*}
is fully faithful.\qed
\end{prop}

\begin{rmk}
A 1-morphism $f:U\to V$ in $\SS$ is faithful (categorically \'etale, a
monomorphism, an 
isomorphism) if and only if the induced 1-morphism $\ul U\to \ul V$ of
presheaves over $\SS$ is faithful (categorically \'etale, a
monomorphism, and isomorphism).  
\end{rmk}

\begin{rmk}
For every presheaf $\FF$ over $\SS$ and every object $U$ of $\SS$,
there is a canonical morphism of groupoids
\begin{equation}\label{choice}
\bhom(\ul U,\FF)\longrightarrow \FF_U\,,
\end{equation}
given by evaluation at $\id_U$.  It is always an equivalence of
groupoids.  Given an object $x$ of $\FF$ lying over $U$, any choice of
pullbacks for $x$ defines a morphism $\ul U\to \FF$ mapping to $x$
under (\ref{choice}).

Thus, it is justified to write $\FF(U)$ instead of $\FF_U$. 
\end{rmk}

\subsubsection{Fibered products of presheaves}

Let $\SS$ be a 2-category.

\begin{prop}
Fibered products exist in the 2-category of presheaves over
$\SS$. The Yoneda functor $\SS\to(\text{\rm presheaves/$\SS$})$ commutes
with any fibered products which exist in $\SS$.
\end{prop}
\begin{pf}
The construction is analogous to the proof of the fact that fibered
products exist in the 2-category of categories fibered in groupoids
over a 1-category.
\end{pf}

\subsection{Topologies and sheaves}

\subsubsection{Sieves}

\begin{defn}
Let $\SS$ be a 2-category and $U$ an object of $\SS$. A {\bf sieve }for $U$
is given by a collection $R$ of objects $V\to U$ of $\ul{U}$ such that
if $V\to U$ is in $R$ and 
$$\xymatrix{
W\rto\drlowertwocell<0>{^<-2>} & V\dto\\
& U}$$
is a 1-morphism in $\ul{U}$, then $W\to U$ is in $R$, also.
\end{defn}

A sieve $R$ for $U$ defines a sub-2-category of $\ul{U}$ by 
$$\bhom_R(W,V)=\bhom_{\ul{U}}(W,V),$$
whenever $V$ and $W$ are objects in $R$.  We shall always identify a sieve
$R$ for $U$ with this sub-2-category of $\ul{U}$ it generates.

Thus, given a sieve $R$ for $\ul{U}$, we have the canonical inclusion
2-functor $R\to \ul{U}$ and by composing with $\ul{U}\to\SS$ a canonical
2-functor $R\to\SS$.

\begin{prop}
If $R$ is a sieve for $\ul{U}$, then $R\to\SS$ is a
presheaf. Conversely, a sub-2-category $R$ of $\ul{U}$, such that $R\to\SS$
is a presheaf comes from a unique sieve for $U$. \qed
\end{prop}

Note that a morphism $V\to U$ is an object of the sieve
$R\subset\ul{U}$ if and only if the induced morphism of
$\SS$-presheaves $\ul{V}\to\ul{U}$ factors through $R\subset\ul{U}$. 

More precisely, if $V\to U$ partakes in the sieve $R\subset\ul{U}$, then
there exists a unique strictly commutative diagram
$$\xymatrix{
{\ul{V}}\rto \drto & R\dto\\
&{\ul{U}}}$$
of presheaves over $\SS$.  If there exists a 2-commutative diagram
$$\xymatrix{
{\ul{V}}\rto \drto \drtwocell\omit{^<-2>}& R\dto\\
&{\ul{U}}}$$
of presheaves over $\SS$, then $V\to U$ partakes in $R$.

\begin{lem}
Finite intersections and arbitrary unions of sieves for $U$ are sieves
for $U$.\qed
\end{lem}

\begin{con}[pullback sieve]
Consider a morphism $f:V\to U$ in $\SS$.  If $R$ is a sieve for $U$,
define $f^{-1}R\subset\ul{V}$ to consist of all $W\to V$ such that the
composition $W\to V\to U$ is in $R$. Of course $f^{-1}R$ is a sieve for
$V$. Note that
$$\xymatrix{
{\phantom{i}f^{-1}R\phantom{i}}\dto \ar@{^{(}->}[r]  & {\ul{V}}\dto\\
{\phantom{i}R\phantom{i}} \ar@{^{(}->}[r] & {\ul{U}}}$$
is a 2-cartesian diagram of presheaves over $\SS$.
\end{con}

\subsubsection{Topologies}

We shall now define topologies on 2-categories.  A topology is
characterized by its collection of {\em covering sieves}. 

\begin{defn}
Let $\SS$ be a 2-category.  A {\bf topology }on $\SS$ is given by the
data

$\bullet$ for every object $U$ of $\SS$ a collection of sieves for
$U$, called the {\em covering sieves }of $U$,

\noindent subject to the constraints:

(i) (pullbacks) for all morphisms $f:V\to U$ in $\SS$ and all covering
sieves $R\subset \ul{U}$ of $U$, the pullback $f^{-1}R\subset\ul{V}$
is a covering sieve for $V$,

(ii) (local nature) if $R\subset\ul{U}$ is a covering sieve for $U$ and
$R'\subset\ul{U}$ is another sieve, which {\em covers $R$-locally},
the sieve $R'$ is also a covering sieve of $U$.  Here we say that $R'$
covers $R$-locally, if for all $f:V\to U$ in
$R$ the pullback $f^{-1}R'$ is a covering sieve of $V$,

(iii) (identities) for every object $U$ of $\SS$ the sieve $\ul{U}$ is
a covering sieve for $U$.

\noindent A 2-category which has been endowed with a topology is
called a {\bf 2-site}.
\end{defn}

Sometimes, a topology can be defined in terms of a pretopology.  A
pretopology is given by its collection of {\em covering families}. 

\begin{defn} \label{defpretop}
Let $\SS$ be a 2-category.  A {\bf pretopology }on $\SS$ is given by
the data

$\bullet$ for every object $U$ of $\SS$ a collection of families of
$U$-objects, called the {\em covering families }of $U$,

\noindent subject to the constraints:

(i) (pullbacks) if $(U_i\to U)_{i\in I}$ is a covering
family of $U$ and $V\to U$ is a morphism, then there exists a
2-pullback $V_i\to V$ of $U_i\to U$, for all $i\in I$, such
that $(V_i\to V)_{i\in I}$ is a covering family of $V$;

(ii) (composition) if $(U_i\to U)_{i\in I}$ is a covering
family of $U$ and for every $i\in I$ we have a covering family
$(V_{ij}\to U_i)_{j\in J_i}$ of $U_i$, then the total family
$(V_{ij}\to U)_{i\in I, j\in J_i}$, obtained by composition, is a
covering family of $U$;

(iii) (identity) the one-member family $(\id:U\to U)$ is a covering
family of $U$, for all objects $U$ of $\SS$.
\end{defn}

\begin{rmk}
Note that in the pullback condition we do not require that {\em every
}pullback of a covering family is a covering family.  We only demand
that there exists at least one pullback family which is covering. 

In particular, a family isomorphic to a covering family need not be
a covering family.
\end{rmk}

We can associate a topology to a given pretopology as follows.  Call a
sieve $R\subset\ul{U}$ {\em covering }if there exists a covering
family (for the given pretopology) $(U_i\to U)$ such that for all
$i$ we have that $U_i$ is in $R$.

\begin{lem}
This defines a topology on $\SS$. This topology is called the topology
{\bf associated }to the given pretopology.
\end{lem}
\begin{pf}
This proof is very similar to the corresponding proof in the
1-category setting.
\end{pf}

\subsubsection{Sheaves}

\begin{defn}
A presheaf $\FF\to\SS$ is called a {\bf sheaf} if, for every object
$U$ of $\SS$ and every covering sieve $R\subset\ul{U}$, the
canonical restriction functor
$$\bhom_{\SS}(\ul{U},\FF)\longrightarrow\bhom_{\SS}(R,\FF)$$
is an equivalence of groupoids.
\end{defn}

Suppose that $(U_i\to U)$ is a covering family for a pretopology. We
denote by $U_{i_0\ldots i_p}$ the fibered product
$U_{i_0}\times_U\ldots\times_U U_{i_p}$.  For a presheaf $\FF$, and an
object $x$ of $\FF_U$, we write $x\resto U_{i_0\ldots i_p}$ for a
chosen pullback of $x$. Given a morphism $\alpha:x\to y$ in $\FF(U)$,
it induces a unique morphism $\alpha\resto U_{i_0\ldots i_p}:x\resto
U_{i_0\ldots i_p}\to y\resto U_{i_0\ldots i_p}$. 

\begin{lem}
Assume that the topology on $\SS$ is defined by a pretopology. Then
the presheaf $\FF\to\SS$ is a sheaf if and only if for every object
$U$ of $\SS$ and every covering family $\uU=(U_i\to U)$ the
following three conditions are satisfied:

(i) Assume given two objects $x,y\in \FF(U)$ and two morphisms
$\alpha,\beta:x\to y$. If $\alpha|U_i=\beta|U_i$, for all $i$, then
$\alpha=\beta$. 

(ii) Assume given two objects $x,y\in \FF(U)$ and for every $i$ a
morphism $\alpha_i:x|U_i\to y|U_i$, such that
$\alpha_i|U_{ij}=\alpha_j|U_{ij}$, for all $i,j$. Then there exists a
morphism $\alpha:x\to y$ such that $\alpha|U_i=\alpha_i$, for all
$i$. 

(iii) Given, for every $i$, an object $x_i$ of $\FF(U_i)$ and for
all $i,j$ a morphism $\alpha_{ij}:x_i\resto U_{ij}\to x_j\resto
U_{ij}$ in $\FF(U_{ij})$, such that for all $i,j,k$ we have
$\alpha_{jk}\resto U_{ijk}\comp\alpha_{ij}\resto U_{ijk}=
\alpha_{ik}\resto U_{ijk}$, then there exists an object $x$ of $\FF(U)$,
and morphisms $\alpha_i:x\resto U_i\to x_i$, such that
$\alpha_{ij}\comp \alpha_i\resto U_{ij}=\alpha_j\resto U_{ij}$, for
all $i,j$. \qed
\end{lem}

\begin{defn}
Given two sheaves $\FF$, $\GG$ over $\SS$, the groupoid of sheaf morphisms
from $\FF$ to $\GG$ is defined to be the groupoid $\bhom_{\SS}(\FF,\GG)$ of
presheaf morphisms from $\FF$ to $\GG$. 
\end{defn}

Thus the sheaves over $\SS$ are a
2-category.  There is the canonical fully faithful inclusion functor 
$$i:(\text{sheaves}/\SS)\longrightarrow(\text{presheaves}/\SS)\,.$$

\begin{prop}\label{sheafmonocrit}
Fibered products of sheaves exist and $i$ commutes with them. A
morphism of sheaves is faithful (categorically \'etale, a
monomorphism, an isomorphism) if 
and only if it is faithful (categorically \'etale, a monomorphism, an
isomorphism) considered 
as a morphism of presheaves (cf.\ Proposition~\ref{monocrit}).\qed
\end{prop}

\begin{defn}
A family $\FF_i\to \FF$ of morphisms of sheaves is called {\bf
epimorphic}, if for every object $U$ of $\SS$ and every morphism $\ul
U\to\FF$, the sieve for $U$, consisting of all $V\to U$ admitting a
2-commutative diagram
$$\xymatrix{
{\ul V}\rto\dto\drtwocell\omit{^} & {\ul U}\dto\\
{\FF_i}\rto & {\FF}}$$
for some $i$, is a covering sieve. 
\end{defn}

\begin{rmk}
If the topology is given by a pretopology, then $\FF_i\to \FF$ is an
epimorphic family of sheaf morphisms, if and only if for every $\ul
U\to\FF$, there exists a covering family $V_j\to U$ for the
pretopology, such that for every $j$ there exists an $i$ and a
2-commutative diagram
$$\xymatrix{
{\ul V_j}\rto\dto\drtwocell\omit{^} & {\ul U}\dto\\
{\FF_i}\rto & {\FF}}$$
\end{rmk}

\subsection{From simplicial categories to 2-categories}

The 2-category we use to construct differential graded schemes comes
from a simplicial category of differential graded algebras.  See
\cite{dgsI} for our conventions concerning simplicial categories and
differential graded algebras.

We may consider every groupoid as a simplicial set, by passing to the
simplicial nerve.  Every groupoid becomes a Kan (i.e.\ fibrant)
simplicial set in this way.  Thus we may consider every 2-category
as a simplicial category.

Recall (see, for example, Page~36 of \cite{simhomthe}) that to a Kan
simplicial set $X$ we 
can associate the fundamental groupoid $\Pi_1 X$ as follows.  Objects of
$\Pi_1 X$ are the vertices $\Delta^0\to X$.  Morphisms in $\Pi_1 X$ are
homotopy classes (relative $\del\Delta^1$) of `paths' $\Delta^1\to X$.
Composition is defined by using the Kan property:  Any two composable
paths give rise to a horn. Filling this horn with a 2-simplex yields
the composition as the third edge.

\begin{prop}
Let $\SS$ be a simplicial category all of whose hom-spaces are
fibrant.  If $U$ and $V$ are objects of $\SS$, define 
$$\bhom(U,V)=\Pi_1\shom(U,V).$$ 
This definition endows $\SS$ with the structure of a 2-category
$\widetilde{\SS}$ in such a way that $\SS\to\widetilde{\SS}$ is a
simplicial functor.\qed
\end{prop}

Note that the condition that $\SS\to\widetilde{\SS}$ is a simplicial
functor determines that structure of 2-category on $\widetilde{\SS}$
uniquely.  It is called the 2-category  {\bf associated }to the
simplicial category $\SS$.

\begin{numex}\label{betc}
Let $\SS$ be a simplicial closed model category and $\UU\subset\SS$ a
full subcateogry all of whose objects are fibrant and cofibrant. Then
all hom-spaces in $\UU$ are fibrant and so there is an associated
2-category $\widetilde{\UU}$.  This is how our
2-categories arise. 
\end{numex}

For the following result concerning the compatability of homotopy
fibered products with fibered products in the associated 2-category,
we need an additional property of the simplicial closed model
category $\SS$. We say that $\SS$ {\em admits finite tensors}, if for
every object $X$ of $\SS$ and every finite simplicial set $K$, there
exists an object $K\otimes X$ in $\SS$, and a simplicial map $K\to
\shom(X, K\otimes X)$, such that for every object $Y$ of $\SS$, the
induced map
$$\Hom(K\otimes X,Y)\longrightarrow\Hom\big(K,\shom(X,Y)\big)$$
is bijective.

\begin{lem} \label{strictistwo}
Let $\SS$ be a simplicial clsoed model category admitting finite
tensors and $\UU$ a full 
subcategory of fibrant-cofibrant objects. Let $V\to U$ be a fibration
in $\UU$, and $U'\to U$ an arbitrary morphism in $\UU$. Consider the
(strict) fibered product in $\SS$
\begin{equation}\label{jcgd}
{\begin{diagram}
V'\dto\rto & U'\dto\\V\rto & U.
\end{diagram}}
\end{equation}
Assume that $V'$ is in $\UU$. Moreover, let at least one of the two
conditions 

(i) for all $Z\in\UU$, we have that $\pi_2\shom(Z,U)=0$,

(ii) for all $Z\in\UU$, we have that
$\pi_1\shom(Z,V')\to\pi_1\shom(Z,U')$ is injective,

be satisfied. Then (\ref{jcgd}) is a 
fibered product in the associated 2-category $\widetilde{\UU}$.
\end{lem}
\begin{pf}
Because $\SS$ admits finite tensors, the diagram
$$\xymatrix{
\shom(Z,V')\dto\rto & \shom(Z,U')\dto\\
\shom(Z,V)\rto & \shom(Z,U)}$$
is a cartesian diagram of simplicial sets, for every object $Z$ of
$\UU$. Moreover, since $V\to U$ is a fibration, by the simplicial
model category axiom, $\shom(Z,V)\to \shom(Z,U)$ is a fibration of
(fibrant) simplicial sets.  Using these two facts, it is easy to prove
that under either of the two assumptions~(i) or~(ii), the induced
diagram
$$\xymatrix{
\Pi_1\shom(Z,V')\dto\rto & \Pi_1\shom(Z,U')\dto\\
\Pi_1\shom(Z,V)\rto & \Pi_1\shom(Z,U)}$$
is a cartesian diagram of groupoids.
\end{pf}

\subsubsection{Resolving algebras}

Recall (see \cite{dgsI}), that a {\em differential graded algebra }is
always graded commutative with unit, over a field $k$ of
characteristic zero. A differential graded algebra is a {\em resolving
algebra }if is free as a graded commutative algebra with unit, on
generators in non-positive degrees.  If finitely many generators in
each degree suffice, we call a resolving algebra {\em quasi-finite},
if in total finitely many generators suffice we speak of a {\em finite
}resolving algebra. A quasi-finite resolving algebra with perfect
complex of differentials is called a {\em perfect }resolving algebra
(see Definition~I.\ref{def.per}). 

We proved (see Corollary~I.\ref{res.fib.cof}), that the resolving
algebras are fibrant-cofibrant 
objects in the simplicial closed model category $\AA$ of all differential
graded algebras.  
Thus, by Example~\ref{betc}, the category of all resolving
algebras admits an associated 2-category, as all
hom-spaces between such algebras are fibrant.

\begin{defn}
This 2-category of resolving algebras is called $\RR$.

The full sub-2-categories of resolving algebras which are
quasi-finite, perfect
or finite are denoted by $\RR_\qf$, $\RR_\perf$  and $\RR_\f$,
respectively. 
\end{defn}

\begin{rmk}
By the results of Section~\ref{sec.homo} in~\cite{dgsI}, a morphism of
resolving algebras $A\to B$ is a quasi-isomorphism if and only if it
is 2-invertible in $\RR$.  If two morphisms $f,g:A\to B$ are
2-isomorphic, they induce identical homomorphisms $h^\ast(A)\to
h^\ast(B)$ on cohomology.
\end{rmk}

\subsubsection{Fibered products}

Since, for the purposes of doing geometry, we will pass to the opposite
category $\RR^\op$ of $\RR$, we will state our results here in terms
of $\RR^\op$. Note that by Remark~I.\ref{fin.tensor}, the opposite of
$\AA$ admits finite tensors.

\begin{prop}\label{abs.prod}
Absolute products exist in the opposite categories of $\RR$, $\RR_\qf$,
$\RR_\perf$ and $\RR_\f$.  The inclusions $\RR_\f\subset
\RR_\perf\subset \RR_\qf\subset \RR$ commute with them. 
\end{prop}
\begin{pf}
Use tensor products over $k$.  Note that $\shom(k,A)=\ast$, for all
differential graded algebras $A$, so that $\pi_2\shom(k,A)=0$, and we
can apply Lemma~\ref{strictistwo}.
\end{pf}

Recall the definitions of {\em \'etale }morphism and {\em standard
\'etale }morphism of quasi-finite resolving algebras,
Definitions~I.\ref{defetale} and~I.\ref{st.et}. 

\begin{prop}\label{base.et.r}
Let $A\to B$ be an \'etale morphism of finite (perfect, quasi-finite)
resolving algebras. Let $A\to A'$ be an arbitrary morphism of finite
(perfect, quasi-finite) resolving algebras.
$$\xymatrix{
& A'\\
B & A\lto_-{\text{\rm\tiny \'etale}}\uto}$$
The induced fibered product in $\RR^\op_\f$ ($\RR^\op_\perf$,
$\RR^\op_\qf$) exists.  If $B'$ is this fibered product, then $A'\to
B'$ is again \'etale. 

Moreover, if $A\to B$ is standard \'etale, then we may choose $A'\to
B'$ to be standard \'etale, too.
\end{prop}
\begin{pf}
By the results of Section~I.\ref{sec.fibered}, we can choose a
finite (quasi-finite) resolution of $A\to B$. Thus, we may
assume without loss of generality that $A\to B$ is itself a
finite (quasi-finite) resolving morphism. We let
$B'=B\otimes_A A'$, which represents a strict fibered product in the
opposite of the category of all resolving algebras.  Moreover, $A'\to
B'$ is again a 	finite (quasi-finite) resolving morphism,
and so $B'$ is a finite (perfect, quasi-finite) resolving
algebra. Clearly, $A'\to B'$ is again \'etale.  By
Proposition~I.\ref{etahom}, Lemma~\ref{strictistwo} applies, and so
$B'$ provides us with a fibered product in $\RR^\op_\f$ ($\RR^\op_\perf$,
$\RR^\op_\qf$).
\end{pf}

\subsection{The \'etale topology}\label{sub.et.top}

We need a base category over which to do geometry.  There are various
choices, all leading to the same notion of differential graded scheme.
This base category will be the opposite category of a suitable category
of resolving algebras, somewhere between $\RR_\f$ and $\RR_\qf$.

\begin{defn}\label{defss}
Let $\SS$  be a full sub-2-category of $\RR^\op$
satisfying

(i) every object of $\SS$ is quasi-finite,

(ii) if $A\to B$ is a finite resolving morphism in $\RR$ and $A$
belongs to $\SS$, then so does $B$.  The ground field $k$ belongs to
$\SS$.
\end{defn}

In particular, all finite resolving algebras are contained in $\SS$.

\begin{numex}
We could let $\SS$ consist of any of the following:

(i) all quasi-finite resolving algebras,

(ii) all perfect resolving algebras,

(iii) all finite resolving algebras.

If we use one of these categories for $\SS$, then we write $\SS_\qf$,
$\SS_\perf$ or $\SS_\f$, respectively.
\end{numex}

We call a morphism $V\to U$ in $\SS$ {\bf \'etale} or {\bf standard
\'etale}, respectively, if the corresponding morphism of quasi-finite
resolving algebras is \'etale or standard \'etale.

We define a functor
\begin{align*}
\SS&\longrightarrow(\text{finite type $k$-schemes})\\*
U&\longmapsto h^0(U)\,,
\end{align*}
by associating to a differential graded algebra $A$ the spectrum of
$h^0(A)$. 
This functor maps \'etale morphisms to \'etale morphisms.

\begin{defn}\label{doetos}
The {\bf \'etale topology }on $\SS$ is defined by
calling, for an object $U$ of $\SS$, a sieve
$R\subset\ul{U}$ {\em covering}, if there exists a family of \'etale
morphisms $U_i\to U$ in $R$ such that $\coprod_i h^0(U_i)\to h^0(U)$
is a surjective morphism of schemes.
\end{defn}

We will show that this notion of covering sieve defines a topology on
$\SS$ by proving that there exists a pretopology  on $\SS$, whose
associated topology is given by Definition~\ref{doetos}.

\begin{defn}\label{doep}
The {\bf \'etale pretopology }on $\SS$ is defined by
calling a family $(U_i\to U)$ a {\em covering family }if

(i) every $U_i\to U$ is standard \'etale,

(ii) $\coprod_i h^0(U_i)\to h^0(U)$ is a surjective morphism of
schemes.
\end{defn}

\begin{prop} 
Definition~\ref{doep} defines a pretopology on $\SS$.
\end{prop}
\begin{pf}
We need to check the three properties of Definition~\ref{defpretop}.

(i) (pullbacks) Let $(U_i\to U)$ be a covering family for the \'etale
pretopology.  Thus every $U_i\to U$ is standard \'etale.  By
Proposition~\ref{base.et.r}, the base change $V_i\to V$ exists in $\SS$
and may be chosen to be standard \'etale, again. Note that $h^0$
commutes with pullback. Hence $(V_i\to V)$ is a covering family for
the \'etale pretopology.

(ii) (composition) This property is satisfied because a composition of
standard \'etale morphisms is standard \'etale.

(iii) the identity property is trivially verified.
\end{pf}

\begin{lem}\label{mailems}
Let $X\to Y$ be an \'etale morphism in $\SS$.  Then there
exists a family of 2-commutative diagrams in $\SS$
\begin{equation}\label{trifact}
{\begin{diagram}
X_i\dto\druppertwocell<0>{<2>} & \\
X\rto & Y
\end{diagram}}
\end{equation}
such that

(i) Every $X_i\to X$ is an open immersion,

(ii) $\coprod h^0(X_i)\to h^0(X)$ is onto,

(iii) every $X_i\to Y$ is standard \'etale.
\end{lem}
\begin{pf}
Translating the Main Lemma~I.\ref{mailem} into the opposite category
$\SS$, we get diagrams
$$\begin{diagram}
X_i'\dto\rto & X_i\dto\\
X\rto & Y
\end{diagram}$$
with 

(i) every $X_i'\to X$ is an elementary open immersion,

(ii) $\coprod h^0(X'_i)\to h^0(X)$ is onto,

(iii) every $X_i'\to X_i$ is an isomorphism in $\RR$,

(iv) every $X_i\to Y$ is standard \'etale.

\noindent Choosing a 2-inverse for $X_i'\to X_i$ we
obtain~(\ref{trifact}). 
\end{pf}

\begin{them}
A sieve $R\subset\ul{U}$ is a covering sieve for the topology induced
by the \'etale pretopology if and only if it satisfies the condition
of Definition~\ref{doetos}.
\end{them}
\begin{pf}
If $R\subset\ul{U}$ is a covering sieve for the associated topology it
satisfies Definition~\ref{doetos} trivially. Let us prove the
converse. Thus assume that $R\subset\ul{U}$ is a sieve and that
$(U_i\to U)_{i\in I}$ is a family of \'etale morphisms in $R$ such
that $\coprod_i h^0(U_i)\to h^0(U)$ is surjective. We have to show
that $R$ contains a covering family for the \'etale pretopology. 

For given $i\in I$, choose a family of 2-commutative diagrams
\begin{equation*}
{\begin{diagram}
V_{ij}\dto\druppertwocell<0>{<2>} & \\
U_i\rto & U,
\end{diagram}}
\end{equation*}
for $j\in J_i$, as in Lemma~\ref{mailems}. Then the total family
$(V_{ij}\to U)_{i\in I,j\in J_i}$ is in $R$ and is a covering family
for the \'etale pretopology.
\end{pf}

\begin{cor}
Definition~\ref{doetos} defines a topology on $\SS$.  
\end{cor}

\begin{cor}
The \'etale topology is the topology associated to the \'etale
pretopology.
\end{cor}

We end this section with a definition:

\begin{defn}
A {\bf differential graded sheaf} is a sheaf on $\SS$ with the \'etale
topology. 
\end{defn}

\Section{Descent theory}\label{sec.desc}

\subsection{Descent for morphisms}

Let us fix a morphism of quasi-finite resolving algebras $C\to A$ and
a family of \'etale 
quasi-finite resolving morphisms $A\to A_i$, such that $\coprod_i\spec
h^0(A_i)\to\spec h^0(A)$ is surjective.  We think of the family
$A\to A_i$ as giving rise to a covering family $U_i\to U$ for the
\'etale topology on $\SS$.

Let $$A_{i_0\ldots i_p}=A_{i_0}\otimes_A\ldots\otimes_A A_{i_p}\,$$
and denote the corresponding object of $\SS$ by $U_{i_0\ldots i_p}$. 
For every  $\lambda:\{0,\ldots,q\}\to\{0,\ldots,p\}$ we have a
canonical morphism of differential graded algebras
\begin{align}\label{lamin}
A_{i_{\lambda(0)}\ldots i_{\lambda(q)}}&\longrightarrow
A_{i_0\ldots i_p}\\* \nonumber
a_0\otimes\ldots\otimes a_q&\longmapsto
\bigotimes_{\lambda(\kappa)=0}a_\kappa\otimes \ldots\otimes
\bigotimes_{\lambda(\kappa)=p}a_\kappa\,. 
\end{align}

Let $C\to B$ be quasi-finite resolving morphism and let
$\sigma=(\sigma_{i_0\ldots i_q})_{i_0,\ldots,i_q}$ be a family of 
$\ell$-simplices  $\sigma_{i_0\ldots i_q}\in\shom_C(B,A_{i_0\ldots
i_q})$.
We denote by 
$$\sigma_{i_{\lambda(0)}\ldots i_{\lambda(q)}}\resto
U_{i_0\ldots i_p}$$
the image of $\sigma_{i_{\lambda(0)}\ldots i_{\lambda(q)}}$ under the
map
$$\shom_C(B,A_{i_{\lambda(0)}\ldots i_{\lambda(q)}})\longrightarrow
\shom_C(B,A_{i_0\ldots i_p})$$
induced by~(\ref{lamin}).

For a composition of maps $\lambda:\{0,\ldots,q\}\to\{0,\ldots,p\}$
and $\mu:\{0,\ldots,p\}\to\{0,\ldots,r\}$ we have
$$\big(\sigma_{i_{\mu\lambda(0)}\ldots i_{\mu\lambda(q)}}\resto
U_{i_{\mu(0)}\ldots i_{\mu(p)}}\big)\resto U_{i_0\ldots i_r} =
\sigma_{i_{\mu\lambda(0)}\ldots i_{\mu\lambda(q)}}\resto
U_{i_{0}\ldots i_{r}}\,.$$
This means that we have a cosimplicial space
$$\xymatrix{
{\displaystyle\prod_{\phantom{j}i\phantom{j}}
\shom_C(B,A_i)}\ar@<.5ex>[r]\ar@<1.5ex>[r] & 
{\displaystyle\prod_{i,j}\shom_C(B,A_{ij})}\ar@<1ex>[r] \ar@<0ex>[r]
\ar@<2ex>[r] & \ldots}$$
For every $\ell\geq0$, we also get a cosimplicial set
\begin{equation}\label{coset}
\xymatrix{
{\displaystyle\prod_{\phantom{j}i\phantom{j}}
\pi_\ell\shom_C(B,A_i)}\ar@<.5ex>[r]\ar@<1.5ex>[r]  
&  {\displaystyle\prod_{i,j}\pi_\ell\shom_C(B,A_{ij})}\ar@<1ex>[r]
\ar@<0ex>[r]  \ar@<2ex>[r] & \ldots\,,}
\end{equation}
where for $\ell\geq1$ this assumes that we have chosen a base point
$P:B\to A$. 

We set
$$H^0\big(\UU,\pi_\ell(B/C)\big)=\ker\big(\lower.3ex\hbox{
\xymatrix{
{\displaystyle\prod_{\phantom{j}i\phantom{j}}
\pi_\ell\shom_C(B,A_i)}\ar@<.5ex>[r]\ar@<1.5ex>[r]
&  {\displaystyle\prod_{i,j}\pi_\ell\shom_C(B,A_{ij})}}}\big)\,.$$
For $\ell\geq1$, having chosen a base point $P:B\to A$, we also define
the pointed set
$$H^1\big(\UU,\pi_\ell(B/C)\big)\,$$ 
in analogy to non-abelian first \v{C}ech cohomology. More precisely, we
set
\begin{multline*}
Z^1\big(\UU,\pi_\ell(B/C)\big)=\\
\{(\alpha_{ij})\in\prod_{ij}\pi_\ell\shom_C(B,A_{ij})\st \forall
i,j,k: \alpha_{ik}\resto U_{ijk}=\alpha_{jk}\resto
U_{ijk}\ast\alpha_{ij}\resto U_{ijk}\}\,.
\end{multline*}
Then we let 
$$C^0\big(\UU,\pi_\ell(B/C)\big)=\prod_i\pi_\ell\shom_C(B,A_i)$$
act (from the left) on $Z^1$ by $(\gamma_i)_i\ast
(\alpha_{ij})_{ij}=(\gamma_j\alpha_{ij}\gamma_{i}^{-1})_{ij}$
and let $H^1\big(\UU,\pi_\ell(B/C)\big)$ by the quotient of $Z^1$ by
this action. 

Finally, for $\ell\geq2$, we associate to (\ref{coset}) the cochain
complex obtained by setting the coboundary map equal to
$\del=\sum_i(-1)^i\del_i$.  We denote the associated cohomology groups
by 
$$H^i\big(\UU,\pi_\ell(B/C)\big)\,.$$

More notation: $V=\spec h^0(A)$, $V_i=\spec h^0(A_i)$ and $\VV$
denotes the
\'etale covering family $V_i\to V$ of affine schemes. Let
$$V_{i_0\ldots i_p}=V_{i_0}\times_V\ldots\times_V V_{i_p}\,.$$
Note that $h^0(U_{i_0\ldots i_p})=V_{i_0\ldots i_p}$. 

\begin{them}[Descent]\label{Descent}
Assume that $B$ is finite over $C$. Then

(i) $H^0\big(\UU,\pi_\ell(B/C)\big)=\pi_\ell\shom_C(B,A)$, for every
$\ell\geq0$, 

(ii) $H^1\big(\UU,\pi_\ell(B/C)\big)=0$, for every $\ell\geq1$,

(iii)  $H^i\big(\UU,\pi_\ell(B/C)\big)=0$, for  all $i\geq2$ and for
every $\ell\geq2$. 
\end{them}
\begin{pf}
Induction on the number $n$ of elements in  a basis for $B$ over $C$. 
Choose a subalgebra $B'\subset B$ such that $B'$ has a $C$-basis of
$n-1$ elements and $B$ has a $B'$-basis consisting of one element $x$
of degree $r$.  By induction, we can assume the theorem to hold for
$B'$. 

Let us start by considering the case $\ell=0$.  Recall
(Section~I.\ref{wcwsa}), that for a given
morphism of $C$-algebras $P:B'\to A$ we have defined the homomorphism
of $h^0(A)$-modules
\begin{align*}
\xi_P:h^{-1}\Deru_C(B',A)&\longrightarrow h^r(A)\\*
D&\longmapsto D(dx)\,.
\end{align*}
The cokernel $\cok\xi_P$ is hence a finitely generated
$h^0(A)$-module.  Thus we may consider $\cok\xi_P$ as a coherent sheaf
on the affine scheme $V$.

Note that the cokernel of 
\begin{align*}
h^{-1}\Deru_C(B',A_{i_0\ldots i_p})&\longrightarrow
h^r(A_{i_0\ldots i_p})\\*
D&\longmapsto D(dx)
\end{align*}
is equal to $\cok\xi_P\otimes_{h^0(A)} h^0(A_{i_0\ldots i_p})$, by
Corollary~I.\ref{prep.desc}. 

Consider the commutative diagram
\begin{equation}\label{ce.di}
\vcenter{\xymatrix{
0\rto & \cok\xi_P \rto\dto & \pi_0\shom_C(B,A)\rto\dto &
\pi_0\shom_C(B',A)\dto \\
0\rto & {C^0(\VV,\cok\xi_P)} \rto & {\prod_i\pi_0\shom_C(B,A_i)}\rto &
{\prod_i\pi_0\shom_C(B',A_i)}}}
\end{equation}
The rows are exact, by Corollary~I.\ref{sim.tr}. The left vertical arrow
is injective by usual \'etale descent theory for coherent modules.
The right vertical arrow is injective by induction hypothesis. 

To  prove injectivity of
$\pi_0\shom_C(B,A)\to\prod_i\pi_0\shom_C(B,A_i)$, we may choose a base
point of $\pi_0\shom_C(B,A)$ and thus a base point $P:B'\to A$, as
above. Thus we have diagram~(\ref{ce.di})  at our disposal, and a
simple chase around the diagram proves the required injectivity.

Now let us  prove surjectivity of 
\begin{equation}\label{sur}
\pi_0\shom_C(B,A)\longrightarrow
\ker\big(\lower.3ex\hbox{
\xymatrix{
{\displaystyle\prod_{\phantom{j}i\phantom{j}}
\pi_0\shom_C(B,A_i)}\ar@<.5ex>[r]\ar@<1.5ex>[r]
&  {\displaystyle\prod_{i,j}\pi_0\shom_C(B,A_{ij})}}}\big) 
\end{equation}
For this we start with the diagram
\begin{equation}\label{l}
\vcenter{\xymatrix{
\pi_0\shom_C(B,A) \rto\dto &\pi_0\shom_C(B',A) \rto\dto &
h^{r+1}(A)\dto \\ 
{\prod_i\pi_0\shom_C(B,A_i)}
\ar@<.5ex>[d]\ar@<-.5ex>[d] \rto
&{\prod_i\pi_0\shom_C(B',A_i)} \ar@<.5ex>[d]\ar@<-.5ex>[d]\rto
&C^0\big(\VV,h^{r+1}(A)\big) \\
{\prod_{ij}\pi_0\shom_C(B,A_{ij})} \rto &
{\prod_{ij}\pi_0\shom_C(B',A_{ij})}  &
}}
\end{equation}
whose right horizontal arrows are defined by evaluation at $dx$.
The rows are exact in the middle by Proposition~I.\ref{ext}. The middle
column is exact by the induction hypothesis and the rightmost column
is injective by usual \'etale descent theory applied to the coherent
sheaf $h^{r+1}(A)$ on $V$.

Let $(\alpha_i)\in\prod_i\pi_0\shom_C(B,A_i)$, such that
$\alpha_i\resto U_{ij}=\alpha_j\resto U_{ij}$, for all $i,j$. Chasing
$(\alpha_i)$ around Diagram~(\ref{l}), we obtain
$\alpha\in\pi_0\shom_C(B,A)$, such that $(\alpha_i)$ and
$(\alpha\resto U_i)$ map to the same element of
$\prod_i\pi_0\shom_C(B',A)$. We also obtain a base point $P:B'\to A$
(to which $\alpha$ maps), supplying us with the diagram
\begin{equation}\label{ce.di.ii}
\vcenter{\xymatrix@C=1.5pc{
0\rto & \cok\xi_P \rto\dto & \pi_0\shom_C(B,A)\rto\dto &
\pi_0\shom_C(B',A)\dto \\
0\rto & {C^0(\VV,\cok\xi_P)} \rto\dto &
{\prod_i\pi_0\shom_C(B,A_i)}\ar@<.5ex>[d]\ar@<-.5ex>[d] \rto &
{\prod_i\pi_0\shom_C(B',A_i)}\\
0\rto & {C^1(\VV,\cok\xi_P)} \rto &
{\prod_{ij}\pi_0\shom_C(B,A_{ij})} & 
}}
\end{equation}
Again, the rows are exact by Corollary~I.\ref{sim.tr}. The leftmost
column is exact by \'etale descent for the coherent sheaf $\cok\xi_P$
on $V$.  Now chasing $(\alpha_i)$ and $(\alpha\resto U_i)$ around
Diagram~\ref{ce.di.ii}, we obtain $\beta\in\cok\xi_P$ such that
$\beta\ast\alpha\resto U_i=\alpha_i$, for all $i$. This proves
surjectivity of (\ref{sur}) and finishes the proof of the theorem in
the case $\ell=0$.

Let us now consider the case $\ell=1$ and prove that
$H^1\big(\UU,\pi_1(B/C)\big) =0$. The advantage over the previous case
is that we now have a fixed base point $P:B\to A$ for all spaces we
consider. We have a commutative diagram
\begin{multline}\label{the.big.nasty}
\vcenter{\xymatrix{
0\rto & {C^0(\VV,\cok\delta)} \rto\dto &
{\prod_i\pi_1\shom_C(B,A_i)}\ar@<.5ex>[d]\ar@<-.5ex>[d]\rto &  \\
0\rto & {C^1(\VV,\cok\delta)} \rto\dto &
{\prod_{ij}\pi_1\shom_C(B,A_{ij})}\ar@<0ex>[d]\ar@<1ex>[d]\ar@<-1ex>[d]
\rto & \\
0\rto & {C^2(\VV,\cok\delta)} \rto &
{\prod_{ijk}\pi_1\shom_C(B,A_{ijk})}\rto & 
}}\\
\vcenter{\xymatrix{
&\pi_1\shom_C(B',A)\dto\rto & \im \delta\rto\dto & 0\\
\rto &{\prod_i \pi_1\shom_C(B',A_i)}\ar@<.5ex>[d]\ar@<-.5ex>[d] \rto
& C^0(\VV,\im\delta)\dto\rto & 0 \\ 
\rto&{\prod_{ij}\pi_1\shom_C(B',A_{ij})}
\rto\ar@<0ex>[d]\ar@<1ex>[d]\ar@<-1ex>[d] 
&C^1(\VV,\im\delta)\rto & 0 \\
\rto &{\prod_{ijk}\pi_1\shom_C(B',A_{ijk})}  &
}}
\end{multline}
Here $\cok\delta$ is the cokernel of the boundary map 
$$\delta:h^{-2}\Deru_C(B',A)\longrightarrow h^{-1}\Deru_{B'}(B,A)=
h^{r-1}(A)\,,$$ 
which is given by $\delta(D)=-D(dx)$ (see
Proposition~I.\ref{boundary}). Moreover, $\im\delta$ is the image of the
boundary map 
$$\delta:h^{-1}\Deru_C(B',A)\longrightarrow h^{0}\Deru_{B'}(B,A)=
h^{r}(A)\,,$$ 
which is given by $\delta(D)=D(dx)$. Both $\cok\delta$ and $\im\delta$ are
coherent sheaves on $V$ and so the first and last columns
of~(\ref{the.big.nasty}) are exact. The $B'$-column is exact by
induction hypothesis.
The rows of (\ref{the.big.nasty}) are exact by Lemma~I.\ref{Main} and
Theorem~I.\ref{Bij}.  We can now prove that the  $B$-column is exact in
the middle by a diagram chase around~(\ref{the.big.nasty}). 

For all other cases of the theorem, note that
$$H^i\big(\UU,\pi_\ell(B/C)\big)=
H^i\big(\VV,h^{-\ell}\Deru_C(B,A)\big)\,,$$  
by Theorem~I.\ref{Bij} and Corollary~I.\ref{prep.desc}.  Thus we are
reduced to usual \'etale descent for the coherent sheaf
$h^{-\ell}\Deru_C(B,A)$ over $V$.
\end{pf}

\begin{cor}
The same holds if we assume only that $B$ is perfect over $C$.
\end{cor}
\begin{pf}
This follows immediately by passing to the limit over the various
truncations $B_{(n)}$.  This is permitted, because of
Corollary~I.\ref{greatc} and also Equation~(I.\ref{limone}), which
features in its proof.  We should remark that for $\ell=0$, we are
only claiming the left exactness of a certain sequence, which is
preserved by taking limits.  

For $\ell=1$, we wish to see that the sequence
$$\xymatrix@C=1.5pc{
0\rto & \pi_1\shom_C(B,A)\rto & C^0\big(\UU,\pi_1(B/C)\big)\ar@{o}[r] &
Z^1\big( \UU,\pi_1(B/C)\big)\rto &0}$$
is exact, in the sense that the group in the middle acts transitively
on the pointed set on the right, in such a way that the stabilizer of
the distinguished point is the group on the left.  This exactness
follows from
$${\projectlim_n}^1\pi_1\shom_C(B_{(n)},A)=0\,,$$
which is true, by Equation~(I.\ref{limone}).
\end{pf}

Notation: if $x,y\in X$ are points of a fibrant simplicial set $X$ and
$\alpha,\beta:x\to y$ are paths in $X$, then we write
$\alpha\sim\beta$ if there exists a homotopy between $\alpha$ and
$\beta$, which fixes the endpoints $x$ and $y$.  In other words,
$\alpha\sim\beta$ if and only if $\alpha$ and $\beta$ define the same
arrow inside the fundamental groupoid $\Pi_1 X$.

\begin{cor}\label{sheaf.prop}
Let $B$ be a perfect resolving algebra.

(i) Given two points $x,y\in\shom(B,A)$ and two paths
$\alpha,\beta:x\to y$ in $\shom(B,A)$, such that for every $i$, we
have $\alpha\resto U_i\sim\beta\resto U_i$, then 
$\alpha\sim\beta$.

(ii) Given two points $x,y\in\shom(B,A)$ and for every $i$ a path
$\alpha_i:x\resto U_i\to y\resto U_i$, such that $\alpha_i\resto
U_{ij}\sim\alpha_j\resto U_{ij}$, for all $i,j$, there exists a path
$\alpha:x\to y$ such that $\alpha\resto U_i\sim\alpha_i$, for all
$i$.

(iii) Given, for every $i$, a point $x_i\in\shom(B,A_i)$, and for all
$i,j$ a path $\alpha_{ij}:x_i\resto U_{ij}\to x_j\resto U_{ij}$, such
that for all $i,j,k$ we have $\alpha_{jk}\resto
U_{ijk}\comp\alpha_{ij}\resto U_{ijk}\sim \alpha_{ik}\resto U_{ijk}$,
there exists a point $x\in\shom(B,A)$ and paths $\alpha_i:x\resto
U_i\to x_i$, such that $\alpha_{ij}\comp \alpha_i\resto U_{ij}\sim
\alpha_j\resto U_{ij}$, for all $i,j$.
\end{cor}
\begin{pf}
This is easy to prove using Theorem~\ref{Descent}.
\end{pf}


\subsection{Hypercubes}

We need a few definitions to make the following more efficient.

\begin{defn}
Let $\SS$ be a 2-category and $I$ a set.  A {\bf truncated
hypercube }in $\SS$ with indexing set $I$ is given by the following
data:

$\bullet$ four families of objects of $\SS$:
$$(U_i)_{i\in I}\quad(U_{ij})_{(i,j)\in I^2}\quad
(U_{ijk})_{(i,j,k)\in I^3}\quad(U_{ijkl})_{(i,j,k,l)\in I^4}$$

$\bullet$ nine families of 1-morphisms as follows:
$$t:U_{ij}\to U_i,\quad s:U_{ij}\to U_j$$
$$p_1:U_{ijk}\to U_{ij},\quad m:U_{ijk}\to U_{ik},\quad p_2:U_{ijk}\to
U_{jk}$$
$$a:U_{ijkl}\to U_{ijk},\quad b:U_{ijkl}\to U_{ijl},\quad 
c:U_{ijkl}\to U_{ikl},\quad d:U_{ijkl}\to U_{jkl}$$

$\bullet$ nine families of 2-isomorphisms as follows:

$-$ three  families of 2-isomorphisms fitting into the truncated cube
\begin{equation}\label{three}
\vcenter{\begin{diagram}
&U_{ijk}\dlto\dto\drto&\\
U_{ij}\dto\drto & U_{ik}\dlto\drto & U_{jk}\dlto\dto\\
U_i&U_j&U_k
\end{diagram}}
\end{equation}

$-$ six families of 2-isomorphisms fitting into half of a hypercube

\begin{equation}\label{six}
\vcenter{\xymatrix{
&&& U_{ijkl}\dllto\dlto\dto\drrto &&\\
& U_{ijk}\dlto\dto\drrto & U_{ijl}\dllto\dto\drrto 
& U_{ikl}\dllto\dlto\drrto && U_{jkl}\dllto\dlto\dto \\
U_{ij} & U_{ik} & U_{il} 
& U_{jk} & U_{jl} & U_{kl}}}
\end{equation}

\noindent subject to the constraint

$\bullet$ all four cubes in the truncated hypercube
$$
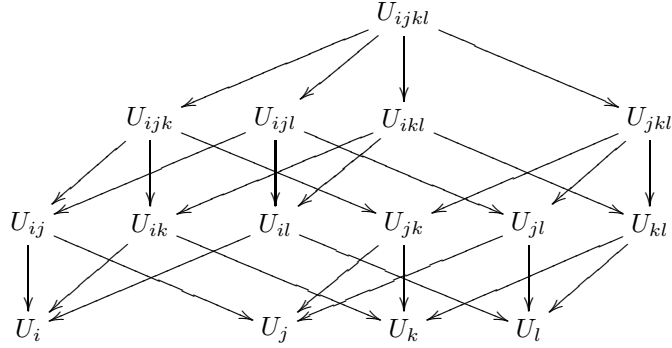
\begin{diagram}
&&& U_{ijkl}\dllto\dlto\dto\drrto &&\\
& U_{ijk}\dlto\dto\drrto & U_{ijl}\dllto\dto\drrto 
& U_{ikl}\dllto\dlto\drrto && U_{jkl}\dllto\dlto\dto \\
U_{ij}\dto\drrto & U_{ik}\dlto\drrto & U_{il}\dllto\drrto 
& U_{jk}\dlto\dto & U_{jl}\dllto\dto & U_{kl}\dllto\dlto \\
U_{i} && U_{j} & U_{k} & U_{l} &
\end{diagram}$$
are 2-commutative, for all $(i,j,k,l)\in I^4$.

We call a truncated hypercube $U\lcom$ {\bf 2-cartesian} ({\bf
strict}), if every one of the 2-commutative squares
appearing in the definition is 2-cartesian (strictly commutative).
\end{defn}

\begin{defn}
Let $U\lcom$ and $V\lcom$ be truncated hypercubes in $\SS$, both with
indexing set $I$.  Then a {\bf morphism }of truncated hypercubes
$\phi\lcom:U\lcom\to V\lcom$ consists of 

$\bullet$ four families of 1-morphisms in $\SS$:
$$\phi_i:U_i\to V_{i},\quad\phi_{ij}:U_{ij}\to
V_{ij},\quad\phi_{ijk}:U_{ijk}\to
V_{ijk},\quad \phi_{ijkl}:U_{ijkl}\to V_{ijkl}\,,$$

$\bullet$ nine families of 2-morphisms fitting into the diagrams:
\begin{equation}\label{mortwo}
\vcenter{\xymatrix@=1pc{
& U_{ij}\ar[dl]\ar[dr]\ar[drrr]^{\phi_{ij}} &&&&\\
U_i\ar[drrr]_{\phi_i} && U_j\ar[drrr]_(.7){\phi_j} &&
V_{ij}\ar[dl]\ar[dr] &\\
&&& V_{i} && V_{j}}}
\end{equation}
\begin{equation}\label{morthree}
\vcenter{\xymatrix{
& U_{ijk}\ar[dl]\ar[d]\ar[dr]\ar[drrr]^{\phi_{ijk}} &&&&\\
U_{ij}\ar[drrr]_{\phi_{ij}} & U_{ik}\ar[drrr]& 
U_{jk}\ar[drrr] && V_{ijk}\ar[dl]\ar[d]\ar[dr]&\\
&&& V_{ij} & V_{ik} &
V_{jk},}}
\end{equation}
\begin{equation}\label{morfour}
\vcenter{\xymatrix{
&& U_{ijkl}\ar[dll]\ar[dl]\ar[d]\ar[dr]\ar[drrr]^{\phi_{ijkl}} &&&&\\
U_{ijk}\ar[drrr]_{\phi_{ijk}}&U_{ijl}\ar[drrr] & U_{ikl}\ar[drrr]& 
U_{jkl}\ar[drrr] &&
V_{ijkl}\ar[dll]\ar[dl]\ar[d]\ar[dr]&\\ 
&&&V_{ijk} & V_{ijl} & V_{ikl} &
V_{jkl},}}
\end{equation}

\noindent subject to the condition that the obvious nine cubes, built
over the nine squares~(\ref{three}) and~(\ref{six}), 2-commute, for
every $(i,j,k)\in I^3$ and  every $(i,j,k,l)\in I^4$.

A morphism of truncated hypercubes is {\bf 2-cartesian} ({\bf
strict}), if every one of the 2-commutative squares
appearing 
in~(\ref{mortwo}), (\ref{morthree})  and~(\ref{morfour})  is
2-cartesian (strictly commutative). 
\end{defn}

\begin{defn}
Given two morphisms $\phi\lcom$ and $\psi\lcom$ from $U\lcom$ to
$V\lcom$, then a {\bf 2-morphism }of morphisms of truncated hypercubes
$\theta\lcom:\phi\lcom\Rightarrow\psi\lcom$ is given by four families
of 2-morphisms in $\SS$: 
$$\xymatrix{U_i\rtwocell^{\phi_i}_{\psi_i}{\phantom{i}\theta_i}
&V_i}\quad  
\xymatrix@C=6ex{U_{ij}\rtwocell^{\phi_{ij}}_{\psi_{ij}}
{\phantom{ij}\theta_{ij}}  
&V_{ij}}\quad 
\xymatrix@C=8ex{U_{ijk} \rtwocell^{\phi_{ijk}}_{\psi_{ijk}}
{\phantom{ijk}\theta_{ijk}} 
&V_{ijk}}\quad 
\xymatrix@C=10ex{U_{ijkl} \rtwocell^{\phi_{ijkl}}_{\psi_{ijkl}}
{\phantom{ijkl}\theta_{ijkl}} &V_{ijkl}}$$
such that the nine families of `2-cylinders' built over~(\ref{mortwo}),
(\ref{morthree}) and~(\ref{morfour}) all 2-commute.
\end{defn}

It is clear that the truncated hypercubes in $\SS$ form a 2-category. 

\begin{defn}
Let $U\lcom$ be a truncated hypercube in $\SS$. An {\bf augmentation
}of $U\lcom$ consists of 

$\bullet$ an object $X$ of $\SS$,

$\bullet$ a family of 1-morphisms:
$$\iota:U_i\longrightarrow X\,$$

$\bullet$ a family of 2-morphisms:
\begin{equation}\label{aug}
\vcenter{
\xymatrix@=1pc{ & U_{ij} \ar[dr]\ar[dl]\ddtwocell\omit{^} & \\
U_i\drto && U_j\dlto\\
& X &}}
\end{equation}

\noindent subject to the constraint that the cube
$$
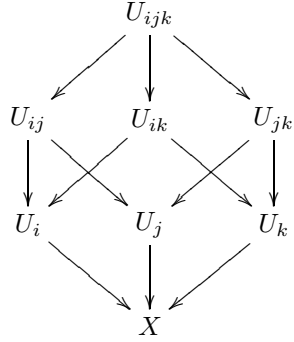
\begin{diagram}
&U_{ijk}\dlto\dto\drto&\\
U_{ij}\dto\drto & U_{ik}\dlto\drto & U_{jk}\dlto\dto\\
U_i\drto &U_j\dto &U_k\dlto\\&X&
\end{diagram}$$
2-commutes, for all $(i,j,k)\in I^3$. 

An augmentation is called {\bf 2-cartesian }({\bf strict}), if 
Diagram~(\ref{aug}) is 2-cartesian (strictly commutative),
for every 
$(i,j)\in I^2$. 

A truncated hypercube endowed with an augmentation is also called a
{\bf hypercube}.  A hypercube is {\bf 2-cartesian }({\bf strict}) if its
underlying truncated hypercube and its augmentation are both 2-cartesian
(strict). 
\end{defn}

\begin{defn}
Let $U\lcom\to V\lcom$ be a morphism of truncated hypercubes. Let
$U\lcom\to X$ and $V\lcom\to Y$ be augmentations.  Then a {\bf morphism
}of augmentations from $X$ to $Y$ consists of 

$\bullet$ a morphism $f:X\to Y$,

$\bullet$ a family of 2-morphisms:
\begin{equation}\label{morau}
\vcenter{
\xymatrix{
U_i\dto\rto^{\phi_i}\drtwocell\omit{^} & V_i\dto \\
X\rto_f & Y}}
\end{equation}

\noindent such that for every $(i,j)\in I^2$ the cube
$$
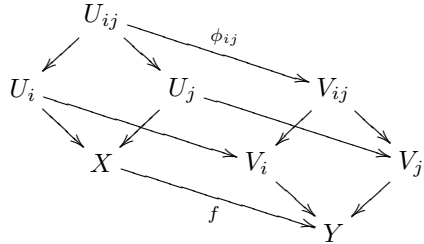
\xymatrix@=1pc{
& U_{ij}\ar[dl]\ar[dr]\ar[drrr]^{\phi_{ij}} &&&&\\
U_i\drto\ar[drrr] && U_j\dlto\ar[drrr] &&
V_{ij}\ar[dl]\ar[dr] &\\
&X\ar[drrr]_f&& V_{i}\drto && V_{j}\dlto\\
&&&& Y &}$$
2-commutes.

A morphism of augmentations is {\bf 2-cartesian }({\bf strict}), if
Diagram~(\ref{morau}) is 2-cartesian (strictly commutative) for every 
$(i,j)\in I^2$. 

A {\bf morphism} of hypercubes is a morphism of truncated hypercubes
together with a morphism of augmentations. A morphism of hypercubes is
{\bf 2-cartesian }({\bf strict}) if both its
underlying morphism of 
hypercubes and its underlying morphism of augmentations are 2-cartesian
(strict). 
\end{defn}

\begin{defn}
Let $(\phi\lcom,f)$ and  $(\psi\lcom, g)$ be two morphisms of
hypercubes.  Then a {\bf 2-isomorphism }of morphisms of augmentations
is a 2-isomorphism
$$\xymatrix{X\rtwocell^f_g{\eta} & Y}\,,$$
such that for every $i\in I$ the `2-cylinder'
$$\xymatrix{
U_i\dto\drrtwocell_{\psi_i}^{\phi_i}{\theta_i} &&\\
X\drrtwocell_g^f{\eta} && V_i\dto\\
&& Y}$$
2-commutes.
\end{defn}

Of course, the hypercubes in $\SS$ also form a 2-category. Moreover,
all the augmentations of a fixed truncated hypercube form a
2-category.


\subsection{Descent for algebras}\label{desforalg}

We begin by describing our setup:

Let $A$ be a quasi-finite resolving algebra and $A\to A_i$ a family of
\'etale quasi-finite resolving morphisms such that $\coprod_i\spec
h^0(A_i)\to \spec h^0(A)$ is surjective. Define $A_{ij}=A_i\otimes_A
A_j$, $A_{ijk}=A_i\otimes_A A_j\otimes_A A_k$ and
$A_{ijkl}=A_i\otimes_A A_j\otimes_A A_k\otimes_A A_l$.
Note that the induced strict hypercube $A\to A\lcom$ in
$\RR_\qf^\op$ is 
2-cartesian, because of  Proposition~\ref{base.et.r}.

Now suppose given perfect resolving morphisms $A_i\to B_i$, $A_{ij}\to
B_{ij}$, $A_{ijk}\to B_{ijk}$ and $A_{ijkl}\to B_{ijkl}$.  Moreover,
assume that these 1-morphisms form part of the data for a strict,
2-cartesian ({\em sic!\/}) morphism of truncated hypercubes $A\lcom\to
B\lcom$ in 
$\RR_\qf^\op$. Finally, assume that the truncated hypercube $B\lcom$
itself is strict (it is necessarily 2-cartesian).

In this situation, we wish to construct 

(i) a perfect resolving morphism $A\to B$,

(ii) for every $i$, a 
morphisms of $A$-algebras 
$f_i:B\to B_i$,

(iii) for all $i,j$, a morphism of $A$-algebras
$f_{ij}:B\to B_{ij}\otimes \Omega_1$, such that 
$$f_{ij}(0)=f_i,\quad f_{ij}(1)=f_j,$$

\noindent in such a way that the following conditions are met:

(iv) for all $i,j,k$ there exists a morphism of $A$-algebras
$f_{ijk}:B\to B_{ijk}\otimes\Omega_2,$ such that
$$f_{ijk}(t,0)=f_{ij}(t),\quad f_{ijk}(0,t)=f_{ik}(t),\quad
f_{ijk}(1-t,t)=f_{jk}(t),$$ 
(in other words, the above data (i), (ii) and (iii) give rise to an
augmentation of the truncated hypercube $B\lcom$)

(v) for every $i$, the (strictly commutative)
square 
\begin{equation}\label{carreaq}
\vcenter{\xymatrix{A\dto\rto & A_i\dto\\
B\rto &B_i}}
\end{equation}
is 2-cartesian.

Our first aim is to reformulate this 2-cartesian requirement. For this,
consider the finite type affine $k$-scheme defined by $h^0(A)$.  It is
endowed with an \'etale cover $h^0(A_i)$. 
The induced cartesian cube in the category of finite type $k$-schemes
is given by 
\begin{equation}\label{notrun}
\xymatrix{h^0(A)\rto &
h^0(A_i)\ar@<-.5ex>[r]\ar@<.5ex>[r] &  
h^0(A_{ij})\ar@<-1ex>[r] \ar@<0ex>[r] 
\ar@<1ex>[r] & h^0(A_{ijk})}\,.
\end{equation}
We also have a cartesian truncated cube
\begin{equation}\label{trun}
\xymatrix{
h^0(B_i)\ar@<-.5ex>[r]\ar@<.5ex>[r] &  
h^0(B_{ij})\ar@<-1ex>[r] \ar@<0ex>[r] 
\ar@<1ex>[r] & h^0(B_{ijk})}\,
\end{equation}
and a cartesian morphism of truncated cubes from~(\ref{notrun})
to~(\ref{trun}). 
In other words, (\ref{trun}) is gluing data for a finite type
$h^0(A)$-algebra $R$:
$$\xymatrix{R\rto &
h^0(B_i)\ar@<-.5ex>[r]\ar@<.5ex>[r] &  
h^0(B_{ij})\ar@<-1ex>[r] \ar@<0ex>[r] 
\ar@<1ex>[r] & h^0(B_{ijk})}\,.$$
Let us denote $h^0(B_\ast)$ by $R_\ast$, for every multi-index $\ast$.

Furthermore, for every $n\leq0$, the truncated cube
$$\xymatrix{
h^n(B_i)\ar@<-.5ex>[r]\ar@<.5ex>[r] &  
h^n(B_{ij})\ar@<-1ex>[r] \ar@<0ex>[r] 
\ar@<1ex>[r] & h^n(B_{ijk})}\,$$
is gluing data for a finitely generated $R$-module $M^n$, because it
is in a certain sense cartesian over~(\ref{notrun}):
$$\xymatrix{M^n\rto &
h^n(B_i)\ar@<-.5ex>[r]\ar@<.5ex>[r] &  
h^n(B_{ij})\ar@<-1ex>[r] \ar@<0ex>[r] 
\ar@<1ex>[r] & h^n(B_{ijk})}\,.$$
Of course, $M^0=R$.

Now, given the data (i), (ii) and (iii), or $f\lcom:B\to
B\lcom$, we get an induced 
morphism of $h^0(A)$-algebras $h^0(B)\to R$ and induced homomorphisms of
$h^0(B)$-modules $h^n(B)\to M^n$ such that for every $i$ the diagram
$$\xymatrix{
h^n(B)\dto\drto &\\
M^n\rto& h^n(B_i)}$$
commutes.

Now we can say that (\ref{carreaq}) is 2-cartesian, if and only if
$h^n(B)\to M^n$ is bijective, for all $n\leq0$.  

If we are given data (i), (ii) and (iii), we say that $f\lcom:B\to
B\lcom$ is a {\bf homotopy square}.  If condition (iv) is satisfied
for $f\lcom:B\to B\lcom$, we say that this homotopy square {\bf
defines an augmentation}.

We will build up $B$ by an inductive procedure.  The following
lemma will be useful:

\begin{lem}\label{great.lemma}
Suppose given a homotopy  square $f\lcom:B\to B\lcom$ and  a finite
resolving morphism $B\to\tilde{B}$ together with a homotopy square
$F\lcom:\tilde{B}\to B\lcom$, such that $F\lcom\resto B$ is equal to
$f\lcom$.  Suppose that $\tilde B$ has a $B$-basis consisting of
finitely many elements, all in the same degree $r$.

Assume that $f\lcom$ defines an augmentation. Then we
may replace the $F_{ij}$ by other morphisms $F'_{ij}:\tilde{B}\to
B_{ij}\otimes\Omega_1$, in such a way that the modified homotopy square
$F'\lcom:\tilde B\to B\lcom$, with $F'_{i}=F_{i}$,
still restricts to $f\lcom$, i.e., $F'\lcom\resto B=f\lcom$, but now
also defines an augmentation.
\end{lem}
\begin{pf}
Choose $f_{ijk}$ as in Condition~(iv).

For purposes of abbreviation, let us introduce the notation
$X_\ast=\shom_A(B,B_\ast)$ and $\tilde
X_\ast=\shom_A(\tilde B,B_\ast)$, for every multi-index $\ast$.
Then we have a commutative diagram of spaces
$$\xymatrix{
{\displaystyle\prod_{\phantom{j}i}\tilde X_i}\dto
\ar@<.5ex>[r]\ar@<1.5ex>[r] & 
{\displaystyle\prod_{i,j}\tilde X_{ij}}\dto \ar@<1ex>[r]
\ar@<0ex>[r]\ar@<2ex>[r] & 
{\displaystyle\prod_{i,j,k}\tilde X_{ijk}}\dto
\ar@<.5ex>[r]\ar@<1.5ex>[r]\ar@<-.5ex>[r]\ar@<2.5ex>[r] &
{\displaystyle\prod_{i,j,k,l}\tilde X_{ijkl}}\dto \\
{\displaystyle\prod_{\phantom{j}i} X_i} \ar@<.5ex>[r]\ar@<1.5ex>[r] & 
{\displaystyle\prod_{i,j} X_{ij}} \ar@<1ex>[r]
\ar@<0ex>[r]\ar@<2ex>[r] & 
{\displaystyle\prod_{i,j,k} X_{ijk}}
\ar@<.5ex>[r]\ar@<1.5ex>[r]\ar@<-.5ex>[r]\ar@<2.5ex>[r] &  
{\displaystyle\prod_{i,j,k,l} X_{ijkl}}}$$
All vertical maps in this diagram are fibrations.

The space $ X_i$ has a canonical base point, given by $f_i$.
The space $ X_{ij}$ has two canonical base points, given by
$f_i$ and $f_j$. Use the notation $ X_{ij}^0$ to denote the
space $ X_{ij}$ endowed with the base point $f_i$ and $ X_{ij}^1$
for $ X_{ij}$ endowed with the base point $f_j$. These two base
points are connected by the path $f_{ij}$, which gives an isomorphism
\begin{equation}\label{wgo}
\pi_\ell  X^0_{ij}\longiso \pi_\ell X^1_{ij}\,,
\end{equation}
for all $\ell\geq0$, which we consider to be {\em canonical}.  The
space $ X_{ijk}$ has three base points, giving rise to three
pointed spaces $ X_{ijk}^0$, $ X_{ijk}^1$ and $
X_{ijk}^2$, with base points $f_i$, $f_j$ and $f_k$, respectively. We
have canonical isomorphisms
$$\pi_\ell  X^0_{ijk}=\pi_\ell  X^1_{ijk}=\pi_\ell 
X^2_{ijk}\,,$$
because, by existence of $f_{ijk}$, it is irrelevant which of the
canonical paths we take between the three different base
points. 

By Theorem~I.\ref{Bij} and its Corollary~I.\ref{greatc}, for
$\ell\geq2$, 
the
homotopy groups $\pi_\ell X_\ast^\epsilon$ are finitely generated
$R_\ast$-modules and by Remark~I.\ref{nature.xi}(iii),
all the canonical maps between them 
are $R_\ast$-linear.
Moreover, we have 
$$\pi_\ell X_i\otimes_{R_i}R_{ij}=\pi_\ell X_{ij}^0$$
and
$$\pi_\ell X_j\otimes_{R_j}R_{ij}=\pi_\ell X_{ij}^1\,.$$
Thus, for $\ell\geq2$, the canonical isomorphisms~(\ref{wgo}) define a
gluing datum for 
a finitely generated $R$-module $ N_\ell$, which comes endowed
with homomorphisms of $R$-modules
$$ N_\ell\longrightarrow \pi_\ell X_i\,,$$
inducing isomorphisms $N_\ell\otimes_R R_i\to \pi_\ell X_i$, and which 
make the diagrams
$$\xymatrix{
{ N_\ell}\rto\drto & {\pi_\ell  X_i}\rto & {\pi_\ell
X^0_{ij}}\dto^{\cong} \\
& {\pi_\ell X_j}\rto & {\pi_\ell  X_{ij}^1} &}$$
commute.  We will only use $N_2$, in what follows.

The $F_i$ induce in a similar way various base points for the $\tilde
X_\ast$. We use notation $\tilde X^\epsilon_\ast$ in a way compatible
with $X^\epsilon_\ast$, to denote the induced  pointed spaces.
Let us denote the fiber of the fibration of pointed spaces
$\tilde X^\epsilon_\ast\to X^\epsilon_\ast$ by $Y^\epsilon_\ast$. 

The $F_{ij}$  induce a commutative diagram
\begin{equation}\label{somme}
\vcenter{
\xymatrix@C=0pc{
\pi_1 Y^0_{ijk}\rrto \drto && \pi_1 Y^1_{ijk} \dlto \\
&\pi_1 Y^2_{ijk}&}}
\end{equation}
because these homotopy groups are abelian (see Lemma~I.\ref{Main}) and
the closed path $\eta_{ijk}=F_{ik}^{-1}\ast F_{jk}\ast F_{ij}$
representing the 
obstruction to commutativity of~(\ref{somme}) maps to the boundary of
a 2-simplex in $X_{ijk}$, and hence can be brought into any fiber
$Y_{ijk}^\epsilon$. 

Thus, by gluing, we obtain  another $R$-module $P_1$
which is locally isomorphic to $\pi_1 Y_i$. There is a canonical
homomorphism of $R$-modules $N_2\to P_1$, which makes the
diagrams
$$\xymatrix{
{ N_2}\dto\rto & {\pi_2  X_i}\dto^{\delta} \\
P_1\rto & \pi_1 Y_i}$$
commute. 

Of particular importance to us will be the $R$-module
$$Q=\cok( N_2\to P_1)\,.$$

Note that $F_{ij}$, $F_{ik}$ and $F_{jk}$ are paths in the space
$\tilde X_{ijk}$, which fit together so as 
to form the circumference of a triangle and thus give rise to an
element of $\pi_1\tilde X_{ijk}^0$, which we shall denote by
$\eta_{ijk}$. Note that under the fibration $\tilde X_{ijk}\to
X_{ijk}$, the homotopy class $\eta_{ijk}$ maps to zero, because its
image forms the boundary of the 2-simplex $f_{ijk}$. Thus, via our
above identifications, we may think of $\eta_{ijk}$ as an element of
$Q\otimes_R h^0(B_{ijk})$.

Now the key observation is that $\eta_{ijk}$ is a \v Cech 2-cocycle
of $\spec R$
with respect to the \'etale covering $\spec h^0(B_i)$ and with values
in the coherent 
$R$-module $Q$. This can be checked by considering the 1-skeleton of a
tetrahedron defined by the $F_{ij}$ inside the space $\tilde
X_{ijkl}$. 

Since this \v Cech cohomology group vanishes, there exist
$$\theta_{ij} \in Q\otimes_R h^0(B_{ij})= \cok(\pi_2 X_{ij}^0\to \pi_1
Y_{ij}^0)=\ker(
\pi_1\tilde X_{ij}^0\to 
\pi_1 X_{ij}^0)\,,$$
such that 
$$\eta_{ijk}=\theta_{ik}^{-1}\ast\theta_{jk}\ast\theta_{ij}\,$$
in $\pi_1 \tilde X_{ijk}^0$.  We are careful to choose representatives
$\theta_{ij}$ which are contained in the fiber $Y_{ij}^0$. 

Now we define
$$F'_{ij}=F_{ij}\ast\theta_{ij}^{-1} \,.$$
More precisely, we choose $F'_{ij}:\tilde{B}\to B_{ij}\otimes\Omega_1$
in such a way that $F'_{ij}$ is homotopic to $ F_{ij}\ast\theta_{ij}^{-1}$
and such that $F'_{ij}\resto B=f_{ij}$.  This is possible because
$\theta_{ij}$ is contained in the fiber of the fibration $\tilde
X_{ij}\to X_{ij}$. 

Now $F'_{ij}$, $F'_{ik}$ and $F'_{jk}$ again form the circumference of
a triangle in $\tilde X_{ijk}$.  The homotopy class of
this triangle is
$$\eta_{ijk}'=\theta_{ik}\ast\eta_{ijk}\ast\theta_{jk}^{-1}
\ast\theta_{ij}^{-1}\,,$$  
which is zero in $\pi_1\tilde{X}_{ijk}^0$, because the kernel of the
homomorphism from 
$\pi_1\tilde X_{ijk}^0$ to $\pi_1X_{ijk}^0$ is abelian.  Thus we can
find a 2-simplex $F'_{ijk}$ in $\tilde X_{ijk}$, whose boundary consists of 
$F'_{ij}$, $F'_{ik}$ and $F'_{jk}$. 

There is no reason why we should be able to make $F'_{ijk}$ restrict
to $f_{ijk}$. For this, $\eta'_{ijk}$ would have to represent zero in
$\pi_1Y_{ijk}^0$, and not just in $\im(\pi_1 Y_{ijk}^0\to\pi_1 \tilde
X_{ijk}^0)$. 
\end{pf}

\begin{them}\label{desalg}
There exists a perfect resolving morphism $A\to B$, together with the
structure of a 2-cartesian  augmentation $f\lcom:B\to B\lcom$, such
that $A\to B$ becomes a strict, 2-cartesian morphism of augmentations
(in $\RR_\qf^\op$). 
\end{them}
\begin{pf}
As mentioned, we will build up $B$ and $f\lcom$ by an inductive
procedure.

Suppose given an integer $n\geq-1$ and a finite resolving morphism
$A\to B_{(n)}$, together 
with morphisms of $A$-algebras $f_i:B_{(n)}\to B_i$ and
$f_{ij}:B_{(n)}\to B_{ij}\otimes\Omega_1$, such that $f_{ij}(0)=f_i$
and $f_{ij}(1)=f_j$. Suppose also that there exists
$f_{ijk}:B_{(n)}\to B_{ijk}\otimes \Omega_2$, satisfying
Condition~(iv), above.

Moreover, assume that the induced  homomorphism of
$h^0(A)$-modules
$$h^\ell(B_{(n)})\longrightarrow M^\ell\,$$ 
is bijective, for all $\ell>-n$, and surjective for $\ell=-n$.

We will construct a resolving morphism $B_{(n)}\to B_{(n+1)}$ and
extensions $F_i$ and $F_{ij}$ of the $f_i$ and the $f_{ij}$ to
$B_{(n+1)}$ in such a way that this extended homotopy square also
defines an augmentation and
$$h^\ell(B_{(n+1)})\longrightarrow M^\ell\,$$ 
is bijective, for all $\ell>-n-1$, and surjective for $\ell=-n-1$.

For the construction, let us choose $f_{ijk}:B_{(n)}\to
B_{ijk}\otimes\Omega_2$, satisfying Condition~(iv). 

We start by choosing elements $b^\nu\in Z^{-n}(B_{(n)})$ whose classes
$[b^\nu]$ in $h^{-n}(B_{(n)})$ generate the kernel of the epimorphism
of $R$-modules $h^{-n}(B_{(n)})\to M^{-n}$ (if $n=-1$, we do not
choose any $b^\nu$, if $n=0$, we choose elements $b^\nu\in B_{(0)}$
generating the kernel of the morphism $h^0(B_{(0)})$-algebras
$h^0(B_{(0)})\to R$). We also choose generators
$m^\mu$ for the $R$-module $M^{-n-1}$ (if $n=-1$, we take generators for
the $h^0(A)$-algebra $R$). Then we choose $\gamma_i^\mu\in Z^{-n-1}(B_i)$
such that $m^\mu$ maps to $[\gamma_i^\mu]$ under $M^{-n-1}\to
h^{-n-1}(B_i)$.

Let ${B_{(n+1)}}=B_{(n)}[x^\nu,y^\mu]$, where $x^\nu$ and $y^\mu$ are
formal variables in degree $-n-1$.  Set $dx^\nu=b^\nu$ and $dy^\mu=0$.
We construct the $F_i$ and $F_{ij}$ by specifying where they send the 
$x^\nu$ and the $y^\mu$.

Let us start with $x^\nu$. Since $b^\nu$ maps to zero in
$h^{-n}(B_i)$, the image $b^\nu_i=f_i(b^\nu)$ of $b^\nu$ in $B_i$ is a
coboundary. Choose $\beta^\nu_i\in B_i^{-n-1}$, such that
$d\beta^\nu_i=b_i^\nu$.  Similarly, let
$b_{ij}^\nu=f_{ij}(b^\nu)$. Since $d(b_{ij}^\nu-b_i^\nu)=0$ and
$(b_{ij}^\nu-b_i^\nu)(0)=0$, by Lemma~I.\ref{basic.sol}, there exists
$\beta^\nu_{ij}\in (B_{ij}\otimes\Omega_1)^{-n-1}$, such that 
$$d\,\beta_{ij}^\nu=b_{ij}^\nu-b_i^\nu\quad\text{and}\quad
\beta_{ij}^\nu(0)=0\,.$$
Again, using Lemma~I.\ref{basic.sol}, we choose $\beta^\nu_{ijk}\in
(B_{ijk}\otimes\Omega_2)^{-n-1}$, 
such that 
$$d\,\beta_{ijk}^\nu(s,t)=
b_{ijk}^\nu(s,t)-b_{ij}^\nu(s)-b_{ik}^\nu(t)+b_i^\nu$$
and
$$\beta_{ijk}^\nu(t,0)=\beta_{ijk}^\nu(0,t)=0\,,$$
where $b^\nu_{ijk}=f_{ijk}(b^\nu)$. 

Let us prove that
$\beta^\nu_{ij}(1)-\beta^\nu_{ik}(1)+\beta^\nu_{jk}(1)$ is a 
coboundary. For this we consider the expression
$$\delta^\nu_{ijk}(t)=\beta^\nu_{ijk}(1-t,t)+
\beta^\nu_{ij}(1-t)+\beta^\nu_{ik}(t)-\beta^\nu_{jk}(t)
-\beta^\nu_{ij}(1)\,,$$ 
which is an element of $B_{ijk}\otimes\Omega_1$.  We have that
$d\,\delta^\nu_{ijk}=0$ and $\delta^\nu_{ijk}(0)=0$, so using
Lemma~I.\ref{basic.sol} once 
again, we find $\Delta^\nu_{ijk}(t)$, such that
$d\Delta^\nu_{ijk}=\delta^\nu_{ijk}$ and $\Delta^\nu_{ijk}(0)=0$. 
Evaluating $\Delta^\nu_{ijk}$ at $t=1$, reveals that
$\beta^\nu_{ij}(1)-\beta^\nu_{ik}(1)+\beta^\nu_{jk}(1)$ is, indeed, a
coboundary.

Thus we have proved that $[\beta^\nu_j-\beta^\nu_i-\beta^\nu_{ij}(1)]$
defines a \v{C}ech 1-cocycle of the affine scheme $\spec R$ and the
\'etale cover $\spec h^0(B_i)$ with values in the coherent sheaf
$M^{-n-1}$. Since this \v{C}ech cohomology group vanishes, we can
bound the cocycle $[\beta^\nu_j-\beta^\nu_j-\beta^\nu_{ij}(1)]$.
So by changing the $\beta^\nu_i$, we may
assume that $[\beta^\nu_j-\beta^\nu_i-\beta_{ij}^\nu(1)]=0$ in
$h^{-n-1}(B_{ij})$.  Therefore,
there exist $\theta^\nu_{ij}\in B^{-n-2}_{ij}$, such that
$d\theta^\nu_{ij}=\beta^\nu_j-\beta^\nu_i-\beta_{ij}^\nu(1)$, for all
$i,j$.  

We now define
\begin{align*}
F_i(x^\nu)&=\beta^\nu_i
\intertext{and}
F_{ij}(x^\nu)(t)&=
(1-t)\beta^\nu_i+ t\big(\beta^\nu_j- \beta_{ij}^\nu(1)\big)+
\beta_{ij}^\nu(t)+(-1)^n\theta^\nu_{ij}\,dt  
\end{align*}

Let us  deal with $y^\mu$.  In this case, $[\gamma^\mu_j-\gamma^\mu_i]$
is directly seen to be zero in $h^{-n-1}(B_{ij})$, and there is no
need (and no freedom anyway) to change the $\gamma^\mu_i$ to find
$\theta^\mu_{ij}\in B^{-n-2}_{ij}$, such that
$d\theta^\mu_{ij}=\gamma^\mu_j-\gamma^\mu_i$, for all $i,j$. 
In this case we set
$$F_i(y^\mu)=\gamma^\mu_i\quad\text{and}\quad
F_{ij}(y^\mu)=
(1-t)\gamma^\mu_i+t\gamma^\mu_j+(-1)^n\theta^\mu_{ij}\,dt\,.$$

We see that this does extend the morphisms of $A$-algebras $f_i$
and $f_{ij}$ to $B_{(n+1)}=B_{(n)}[x^\nu,y^\mu]$, and that the
relations
$$F_{ij}(0)=F_i\quad\text{and}\quad F_{ij}(1)=F_j$$
are satisfied. Moreover, by construction,
$h^\ell({B_{(n+1)}})\longrightarrow M^\ell$ is 
bijective, for all $\ell>-n-1$, and surjective, for all
$\ell\geq-n-1$. 

The only thing left to worry about is if the homotopy square
$F\lcom:B_{(n+1)}\to B\lcom$ defines an augmentation. If it does not,
then we apply Lemma~\ref{great.lemma}.  This leads to a change in the
$F_{ij}$, but since the $F_i$ are not affected, the properties of 
$h^\ell({B_{(n+1)}})\longrightarrow M^\ell$ are not affected. 

Thus our inductive procedure works.  Starting with $B_{(-1)}=A$, we
let 
$$B=\injectlim_n B_{(n)}\,.$$
Since at every step of the induction $F_i$ and $F_{ij}$ are extensions
of $f_i$ and $f_{ij}$, we get induced morphisms of $A$-algebras
$f_i:B\to B_i$ and $f_{ij}:B\to B_{ij}\otimes \Omega_1$, satisfying
$f_{ij}(0)=f_i$ and $f_{ij}(1)=f_j$.  Therefore, we also get
an induced morphism of $h^0(A)$-algebras $h^0(B)\to R$ and induced 
homomorphisms of $h^0(B)$-modules $h^\ell(B)\to M^\ell$, for all
$\ell\leq0$. All these are isomorphisms, by construction, and hence
the (strictly commutative) square
$$\xymatrix{
A\rto\dto & A_i\dto\\
B\rto & B_i}$$
is 2-cartesian in $\RR_\qf^\op$.
This proves that $A\to B$ is perfect, as this property is local in the
\'etale topology on $A$.

To prove that we have an augmentation $f\lcom:B\to B\lcom$, we need to
check that for all $i,j,k$ the path $f_{ik}^{-1}\ast f_{jk}\ast
f_{ij}$ in $\pi_1\big(\shom_A(B,B_{ijk}),f_i\big)$ represents
zero.  But by construction, the image of this path
in $\pi_1\big(\shom_A(B_{(n)},B_{ijk}),f_i\big)$ represents zero, for all
$n$. Thus we conclude, using the fact that 
$$\pi_1\shom_A(B,B_{ijk})=\projectlim_n
\pi_1\shom_A(B_{(n)},B_{ijk})\,$$
which is Corollary~I.\ref{greatc}.
\end{pf}


\subsection{Strictifying truncated hypercubes}

\begin{defn}
A morphism $f:X\to Y$ in a 2-category $\SS$ is called a {\bf
fibration}, if for every 2-commutative diagram
$$\xymatrix{
&& X\dto^f\\
T\rrto_y\urruppertwocell<3>^x{\omit}\rrtwocell\omit{^<-3>\eta}& & Y}$$
there exists a lift $\eta'$ of $\eta$ to $X$, i.e., a diagram
$$\xymatrix{
&& X\dto^f\\
T\rrto_y\urrtwocell^x_{x'}{^\eta'\phantom{l}}& & Y}$$
such that $f\comp\eta'=\eta$.

Another way to say this is that the morphism of groupoids $X(T)\to
Y(T)$ makes $X(T)$ into a fibered category over $Y(T)$, for all
objects $T$ of $\SS$. 
\end{defn}

\begin{defn}
Let  $\SS$ be a 2-category with a distinguished class of 1-morphisms
called $F$-morphism.  We say that $\SS$ {\bf has enough
}$F$-morphisms, if for
every morphism $X\to Y$ in $\SS$ there exists a strict factorization 
$$\xymatrix{X\drto\rto^{\sim} & X'\ar@{->>}[d] \\
& Y}$$
where $X\to X'$ is 2-invertible and $X'\to Y$ is an $F$-morphism.
We call any such factorization an {\bf $F$-resolution }of $X\to Y$
\end{defn}

\begin{prop}\label{strictify}
Let $\SS$ be a 2-category with a distinguished class of morphisms
called $F$-morphism such that

(i) every $F$-morphism is a fibration,

(ii) strict fibered products exist in $\SS$, if at least one of the
two participating morphisms is of type $F$,

(iii) A strict base change of an $F$-morphism is an $F$-morphism,

(iv) compositions of $F$-morphisms are $F$-morphisms.

\noindent Then if $U\lcom$ is a strict truncated hypercube and $V\lcom\to
U\lcom$ is a morphism of truncated hypercubes, then there exists a
strictly commutative diagram
$$\xymatrix{V\lcom\rto^{\sim}\drto & V\lcom'\dto\\
& U\lcom}$$
where $V\lcom \to V\lcom'$ is a 2-invertible morphism of truncated
hypercubes, and $V\lcom'$, as well as $V\lcom'\to U\lcom$, is strict.
Moreover, $V\lcom'\to U\lcom$ 
can be chosen such that all its structure 1-morphisms are
$F$-morphisms.
\end{prop}
\begin{pf}
Start by choosing $F$-resolutions $V_i\to
V_i'\to U_i$ of $V_i\to U_i$. Then, for every $(i,j)\in I^2$, replace
the two morphisms $V_{ij}\to V_i'$ and $V_{ij}\to V_j'$ by
2-isomorphic ones, in such a way that the two squares 
\begin{equation}\label{dxk}\vcenter{
\xymatrix{V_{ij}\rto\dto & V_{i}'\dto\\
U_{ij}\rto & U_{i}}}\quad\vcenter{\xymatrix{V_{ij}\rto\dto & V_{j}'\dto\\
U_{ij}\rto & U_{j}}}
\end{equation}
commute strictly.

Next, consider the strict fibered products
$${\xymatrix{P_{ij}\dto\rto & {V_i'\times V_j'}\dto\\
U_{ij}\rto & {U_i\times U_j}}}$$
We have canonical morphisms $V_{ij}\to P_{ij}$, which we $F$-resolve:
$V_{ij}\to V_{ij}'\to P_{ij}$. 
We replace $V_{ij}$ by $V_{ij}'$. This preserves the strictness of the
diagrams~(\ref{dxk}). It makes $V_{ij}'\to U_{ij}$ into
$F$-morphisms. 

Use the fact that $V_{ij}'\to U_{ij}$ is a fibration to replace
$V_{ijk}\to V_{ij}'$ by a 2-isomorphic morphism making the square
$$\xymatrix{
V_{ijk}\dto\rto & V_{ij}'\dto\\
U_{ijk}\rto & U_{ij}}$$
strictly commutative.

Next, consider the diagram
\begin{equation}\label{sjkfdl}
\vcenter{\xymatrix{
V_{ijk}\drto\rto\dto\ddrtwocell\omit\drrtwocell\omit{^}  &
V_{ij}'\drto & \\ 
U_{ijk}\drto& V_{ik}'\dto\rto & V_{i}'\dto\\
& U_{ik}\rto & U_i}}
\end{equation}
Both the square and the exterior hexagon in this diagram strictly
commute.  Using the fact that $V_{ik}'\to U_{ik}\times_{U_i} V_i'$ is
a fibration, we may replace $V_{ijk}\to V_{ik}'$ by a 2-isomorphic
morphism making the whole diagram (\ref{sjkfdl}) strictly
commute.

Thus, we now have a diagram
\begin{equation}\label{skfdl}
\vcenter{\xymatrix{
V_{ijk}\ddrtwocell\omit\drrtwocell\omit{^}\drto\dto\rto  &
V_{ij}'\times V_{ik}'\drto & \\ 
U_{ijk}\drto & V_{jk}'\dto\rto & V_{j}'\times V_k'\dto\\
& U_{jk}\rto & U_j\times U_k}}
\end{equation}
in which, again,  both the square and the exterior hexagon strictly
commute. We exploit the fact that $V_{jk}'\to P_{jk}$ is a fibration,
to make (\ref{skfdl}) strictly commute.  At this point all squares in
$V\lcom$, as well as all squares in $V\lcom\to U\lcom$, which have
$V_{ijk}$ as source, are strictly commutative. 

The next step is to consider the strict fibered products
$${\xymatrix{P_{ijk}\dto\rto & {V_{ij}'\times V_{ik}'\times V_{jk}'}\dto\\
U_{ijk}\rto & {U_{ij}\times U_{ik}\times U_{jk}}}}$$
and to resolve the canonical morphism $V_{ijk}\to P_{ijk}$ by the
composition $V_{ijk}\to V_{ijk}'\to P_{ijk}$. As above, we can
strictify all squares in $V\lcom$ and in $V\lcom\to U\lcom$ whose
source is $V_{ijkl}$. 

Finally, resolve $V_{ijkl}\to P_{ijkl}$, do finish the proof.
\end{pf}

We will apply this proposition to $\RR^\op_\perf$, with `$F$-morphism'
meaning `perfect resolving morphism'.


\Section{Differential graded schemes}\label{sec.schemes}

\subsection{Differential graded sheaves}

We return to considering the 2-category $\SS$ with the \'etale
topology, introduced in Section~\ref{sub.et.top}.  Recall that we had
defined a {\em differential graded sheaf }to be any sheaf on $\SS$.

If $B$ is a perfect resolving algebra, then it induces a representable
presheaf $X_\qf$ over $\SS_\qf$.  Let us restrict $X_\qf$ to $\SS$
via the the canonical embedding $\SS\to\SS_\qf$. We get a presheaf
$X$ over $\SS$, which may not be representable (depending on the
choice of $\SS$) if $B$ is not finite.

By Corollary~\ref{sheaf.prop}, this presheaf $X$ over $\SS$ is  
a sheaf.  We denote the differential graded sheaf $X$ by 
$$\dgspec B\,.$$

If $B\to C$ is a morphism of perfect resolving algebras, then we get
an induced morphism of differential graded sheaves 
$$\dgspec C\longrightarrow\dgspec B\,.$$
In fact, we get a contravariant 2-functor
$$\dgspec:\RR_\perf\longrightarrow(\text{differential graded
sheaves})\,.$$

\begin{lem}\label{loc-fin-aff}
Let $B$ be a perfect resolving algebra and $U$ the object of
$\SS_\perf$ given by $B$.  Then there exists a finite
covering family of $U_i\to U$ for the \'etale pretopology on
$\SS_\perf$, such that every $U_i$ is given by a finite resolving
algebra.
\end{lem}
\begin{pf}
By Theorem~I.\ref{loc.fin} there exist $g_i\in B^0$, such that the
elementary open immersions $B\to B_{\{g_i\}}$ define a covering family
and 
such that each $B_{\{g_i\}}$ is quasi-isomorphic to a finite resolving
algebra.
\end{pf}

\begin{cor}\label{comp.l}
In case $\SS$ is contained in $\SS_\perf$, the canonical restriction
2-functor 
$$(\text{\rm sheaves over $\SS_\perf$})\longrightarrow (\text{\rm
sheaves over $\SS$})$$
is fully faithful.
\end{cor}
\begin{pf}
By Lemma~\ref{loc-fin-aff}, $\SS$ generates $\SS_\perf$.
\end{pf}

\begin{cor}\label{dgsff}
The contravariant 2-functor\/ 
$$\dgspec:\RR_\perf\to(\text{\rm differential graded
sheaves})$$  
is fully faithful.
\end{cor}
\begin{pf}
First let us consider the case that $\RR_\perf\subset\SS$. 
By Yoneda's lemma, Proposition~\ref{2Y}, we have a fully
faithful 2-functor
$$\RR_\perf\longrightarrow(\text{presheaves over $\SS$})\,.$$
By Corollary~\ref{sheaf.prop}, this 2-functor maps into the
subcategory of sheaves over $\SS$.  So we have a fully faithful
2-functor
$$\RR_\perf\longrightarrow(\text{sheaves over $\SS$})\,.$$

Now let us consider the case that $\SS\subset\RR_\perf$.  The above
considerations applied to $\SS_\perf$ give us a fully faithful 2-functor
$$\RR_\perf\longrightarrow(\text{sheaves over $\SS_\perf$})\,.$$
Composing with the fully faithful 2-functor of Corollary~\ref{comp.l}
finishes the proof.
\end{pf}

\subsection{Affine differential graded schemes}

\begin{defn}
A differential graded sheaf, which is isomorphic to $\dgspec B$, for
some perfect resolving algebra $B$, is called an {\bf affine
differential graded scheme}, or simply {\bf affine}.

If an affine differential graded scheme is isomorphic to $\dgspec B$,
for a finite resolving algebra $B$, we call it {\bf finite affine}, or
simply {\bf finite}. 
\end{defn}

\begin{numnote}\label{equivsep}
By Corollary~\ref{dgsff}, the 2-functor $\dgspec$ induces 
contravariant equivalences of 2-categories
\begin{align*}
\dgspec:\RR_\perf&\longrightarrow(\text{affine differential graded
schemes}) \\
\dgspec:\RR_\f&\longrightarrow(\text{finite affine differential graded
schemes})\,.
\end{align*}
\end{numnote}

\begin{numex}
Every complete intersection in affine space over $k$ gives rise to an
affine differential graded scheme, which is well-defined up to
isomorphism. This is because every complete intersection in affine
space is defined by a regular sequence, and so the associated Koszul
complex gives a differential graded algebra, which is determined by
the complete intersection scheme up to quasi-isomorphism.
For example, any finite extension field $K$ of $k$ may be considered
as an affine differential graded scheme. 

If $A$ is the usual affine coordinate ring of such a complete
intersection, 
we commit the abuse of writing $\dgspec A$ for this differential
graded scheme. 
\end{numex}

\begin{lem}
If $X$ and $Y$ are affine differential graded schemes, then so is
$X\times Y$.  If $X$ and $Y$ are finite, then so is
$X\times Y$. 
\end{lem}
\begin{pf}
This follows from Proposition~\ref{abs.prod}.
\end{pf}

Let $f:X\to Y$ be a morphism between affine differential graded
schemes, where $X\cong\dgspec B$ and $Y\cong\dgspec C$. Then, by
Corollary~\ref{dgsff},  there exists a morphism of perfect resolving
algebras $\phi:C\to B$, such that $\dgspec\phi=f$. The morphism $\phi$
is unique up to homotopy.

\begin{defn}
The morphism of affine differential graded schemes $f:X\to Y$ is called
{\bf \'etale}, (an {\bf open immersion}) if the corresponding morphism
$\phi:C\to B$ of perfect 
resolving algebras is \'etale (an open immersion).  (See
Definitions~I\ref{defetale} and~I.\ref{def.op.imm}.) 
\end{defn}

\begin{lem}\label{prodet}
If $f:X\to Y$ is a morphism of affine differential graded schemes
and $Y'\to Y$ is an \'etale morphism (an open immersion) of affine
differential graded 
schemes, then the fibered product $X'=X\times_{Y}Y'$ is an
affine differential graded scheme and $X'\to X$ is \'etale (an open
immersion).  If 
$X$, $Y$ and $Y'$ are finite, then so is $X'$. 
\end{lem}
\begin{pf}
This follows directly from
Proposition~\ref{base.et.r}.
\end{pf}

\begin{lem}[affine descent]\label{app-rep}
Let $\XX\to U$ be a morphism of differential graded sheaves, where $U$
is affine. Assume that there exists an epimorphic family of morphisms
of differential graded sheaves $U_i\to U$, where every $U_i$ is
affine, and every morphism $U_i\to U$ is \'etale.  Assume that for
every $i$, the 2-fibered product $X_i=\XX\times_U U_i$ is affine.
Then $\XX$ itself is affine.
\end{lem}
\begin{pf}
Without loss of generality, the morphisms $U_i\to U$ are given by
(\'etale) quasi-finite resolving morphisms $A\to A_i$ of perfect
resolving algebras $A$, $A_i$. Because $(U_i\to U)$ is epimorphic,
$\coprod_i\spec h^0(A_i)\to\spec h^0(A)$ is surjective. Define the
strict cartesian hypercube $A\to A\lcom$ in $\RR_\perf^\op$ as in
Section~\ref{desforalg}. Let $U\lcom=\dgspec A\lcom$, so that
$U\lcom\to U$ is a (strict, cartesian) hypercube of affine
differential graded schemes. 

Form the fibered products $X\lcom=\XX\times_U U\lcom$ of differential
graded sheaves. These form a 2-cartesian morphism
$$X\lcom\to U\lcom$$
of 2-cartesian truncated hypercubes of differential graded sheaves. By
Lemma~\ref{prodet}, the truncated hypercube $X\lcom$ consists of
affine differential graded schemes. Thus, via Corollary~\ref{dgsff},
we may choose a 2-cartesian truncated hypercube $B\lcom$ of perfect
resolving algebras, such that $\dgspec B\lcom=X\lcom$.  The morphism
$X\lcom\to U\lcom$ gives rise to a 2-cartesian morphism $A\lcom\to
B\lcom$. 

Now apply Proposition~\ref{strictify} to $A\lcom\to B\lcom$ in
$\RR^\op_\perf$, with `$F$-morphism' meaning `perfect resolving
morphism', to show that we may assume, without loss of generality,
that 
$B\lcom$, as well as $A\lcom\to B\lcom$, is strict.  Moreover, we may
assume that $A\lcom\to B\lcom$ consists of perfect resolving
morphisms.
Thus we are now in the setup of Section~\ref{desforalg}, and from
Theorem~\ref{desalg} we obtain a perfect resolving morphism $A\to B$
and a 2-cartesian morphism of hypercubes $[A\to A\lcom]\to[B\to
B\lcom]$ in 
$\RR^\op_\perf$. Applying $\dgspec$, we obtain a cartesian morphism of
hypercubes $[X\lcom\to X]\to[U\lcom\to U]$ of affine differential
graded schemes.

Then $X\cong \XX$ as sheaves on
$\SS$, and so $\XX$ is representable by the perfect resolving algebra
$B$. 
\end{pf}

\begin{rmk}
We cannot conclude that $\XX$ is finite, even if $U$, all $U_i$ and
all $X_i$ are finite. A counterexample can be derived from
Example~I.\ref{abstract} or, more specifically,
Example~I.\ref{concrete}.  Let $A$ and $B$ be the 
differential graded algebras defined in Example~I.\ref{concrete}. Thus,
$$A=k[x,y,]\{\xi\}/d\xi=y^2-4(x^3-x)\,,$$ 
in the notation of Section~I.\ref{sec.loc.str}. Moreover, $B$ is a
quasi-finite resolution of the differential graded algebra $h^0(A)\oplus
L$. Here $L$ is a projective $h^0(A)$-module of rank one, which we
put in degree -1.  The differential on $h^0(A)\oplus L$ is zero.  We
get a counterexample if $L$ is non-trivial. 

Suppose the elementary open immersion $A\to A_i$ trivializes $L$. Then
$B_i=B\otimes_A A_i$ is quasi-isomorphic to a finite resolving
algebra.  Thus $A$, all $A_i$ and all $B_i$ are `essentially' finite,
but $B$ is not.
\end{rmk}

\begin{prop}
A morphism of affine differential graded schemes is \'etale if and
only if it is a categorically \'etale morphism of sheaves.
It is an open immersion if and only if it is a monomorphism of
sheaves (cf. Propositions~\ref{monocrit} and~\ref{sheafmonocrit}).
\end{prop}
\begin{pf}
Let our morphism be $\dgspec B'\to \dgspec B$.
First we reduce to the case that $B\to B'$ is a resolving
morphism. Then the 
categorically \'etale property translates into bijectivity of
$$\pi_1\shom(B',A)\longrightarrow\pi_1\shom(B,A)\,,$$
for all resolving morphisms $B'\to A$.  By Proposition~I.\ref{etahom},
this is equivalent to  $B\to B'$ being \'etale.

For the second claim, we may assume that $\dgspec B'\to \dgspec B$ is
\'etale. Then $\dgspec B'\to \dgspec B$ is an open immersion, if and
only if the diagram
$$\xymatrix{
{\dgspec B'}\dto\rto & {\dgspec B'}\dto \\
{\dgspec B'}\rto & {\dgspec B}}$$
is 2-cartesian. This is the case if and only if 
$$\pi_0\shom(B',A)\to\pi_0\shom(B,A)$$
is injective, for all $A$, which is the monomorphism property of
$\dgspec B'\to \dgspec B$.
\end{pf}

\begin{defn}
An {\bf affine Zariski cover }of an affine differential graded scheme
$X$ is a collection of open immersions $U_i\to X$, where every $U_i$
is affine. 
\end{defn}

\begin{prop}\label{an-loc}
Every affine differential graded scheme $X$ admits an affine Zariski
cover
$U_i\to X$, such that for every multi-index $(i_0,\ldots,
i_p)$, with $p\geq0$, the fibered product
$$U_{i_0\ldots i_p}=U_{i_0}\times_X\ldots\times_X
U_{i_p}$$ 
is finite.
\end{prop}
\begin{pf}
Let $X=\dgspec B$. As in Lemma~\ref{loc-fin-aff}, we choose $g_i\in
B^0$, such that the 
elementary open immersion $B\to B_{\{g_i\}}$ cover $X$, and
such that each $B_{\{g_i\}}$ is quasi-isomorphic to a finite resolving
algebra $A_i$. 

The fibered product $\dgspec B_{\{g_i\}}\times_X\dgspec
B_{\{g_j\}}$ is represented by $B_{\{g_ig_j\}}$, which is
quasi-isomorphic to ${A_i}_{\{g_j\}}$ and hence finite.  Similarly
for iterated fibered products.
\end{pf}

\subsection{Affine \'etale morphisms}

\begin{defn}
A morphism $\FF\to\GG$ of differential graded sheaves is called {\bf
affine \'etale} (an {\bf affine open immersion}), if for every
morphism $U\to\GG$, with $U$ affine, the 
fibered product $V=\FF\times_{\GG}U$ is affine and the
morphism $V\to U$ is \'etale (an open immersion).
\end{defn}

\begin{prop}
A morphism $\FF\to \GG$ is affine \'etale (an affine open immersion)
if and only if there exists 
an epimorphic  family $U_i\to\GG$ of morphisms with (finite) affine $U_i$,
such that  for every $i$ the 
fibered product $V_i=\FF\times_\GG U_i$ is affine and the morphism
$V_i\to U_i$ is \'etale (an open immersion).
\end{prop}
\begin{pf}
This follows directly from Lemmas~\ref{prodet} and~\ref{app-rep}.
\end{pf}

In particular, a morphism of affine differential graded
schemes is \'etale if and only if it is affine \'etale. Similarly for
open immersions.

\begin{prop}
If $\FF\to\GG$ and $\FF'\to\GG'$ are affine \'etale morphisms of
differential graded sheaves, then so is $\FF\times \FF'\to\GG\times
\GG'$. 
\end{prop}
\begin{pf}
This follows from the fact the affine \'etale property is stable under
composition and arbitrary base change.
\end{pf}

\begin{numnote}\label{loclet}
Let $f:X\to Y$ be a morphism of affine differential graded
schemes. Suppose there exists an epimorphic family of affine \'etale
morphisms $U_i\to X$, where for all $i$ the composition $U_i\to Y$ is
\'etale.  Then $f$ is \'etale.
\end{numnote}

\subsection{Differential graded schemes}

\begin{defn}
A differential graded sheaf $\XX$ is called a {\bf
differential graded scheme}, if there exists an epimorphic family of
affine \'etale morphisms 
${{U}}_i\to\XX$ such that each $U_i$ is affine.

Any such epimorphic family $U_i\to \XX$ is called an {\bf affine
\'etale cover 
}of $\XX$. If all $U_i\to \XX$ are affine open immersion, we speak of
an {\bf affine Zariski cover }of $\XX$.
\end{defn}

\begin{rmk}
It would be more accurate to call these objects 'differential graded
algebraic spaces with affine diagonal', but we find that terminology
too clumsy. 

It seems likely that one can iterate this definition, and obtain more
general objects, which would be differential graded schemes with
weaker separation condition.  At this point it is not clear how useful
this would be, and so we call the above objects simply differential
graded schemes, without a further qualifier. 
\end{rmk}

\begin{prop}\label{dglfin}
If $\XX$ is a differential graded scheme, then there exists an
affine \'etale cover $U_i\to \XX$, where every $U_i$ is finite.
\end{prop}
\begin{pf}
This follows directly from Proposition~\ref{an-loc}.
\end{pf}

\begin{prop}
Let $\XX$ and $\YY$ be differential graded schemes. Then
$\XX\times\YY$ is a differential graded scheme.
\end{prop}
\begin{pf}
Let $U_i\to\XX$ and $V_j\to\YY$ be affine \'etale covers. Then
$U_i\times V_j\to \XX\times\YY$ is an affine \'etale cover of the
product $\XX\times \YY$.
\end{pf}

\begin{lem}[descent]\label{descjdskfl}
Let $\FF\to\XX$ be a morphism of differential grade sheaves, where
$\XX$ is a differential graded scheme. Let $U_i\to\XX$ be an affine
\'etale cover, such that the fibered product
$\VV_i=\FF\times_{\XX}U_i$, is a differential graded scheme, for all
$i$. Then $\FF$ is a differential graded scheme.
\end{lem}
\begin{pf}
Let $V_{ij}\to \VV_i$ be an affine \'etale cover of $\VV_i$, for all
$i$. Since the affine \'etale property, as well as the epimorphism
property are stable under base change and composition, it follows that
$V_{ij}\to\FF$ is an affine \'etale cover. 
\end{pf}

\begin{defn}\label{def-amp}
A differential graded scheme $\XX$ is of {\bf amplitude }$N$, if for
every affine \'etale morphism $\dgspec B\to \XX$,   the
perfect resolving algebra $B$ is of amplitude $N$. 
\end{defn}

If there exists an affine \'etale cover $\dgspec B_i\to \XX$, where
for every $i$, the perfect resolving algebra $B_i$ is of amplitude
$N$, then $\XX$ is of amplitude $N$.

\subsection{Affine morphisms}

\begin{defn}
A morphism $\XX\to\YY$ of differential graded schemes is called {\bf
affine}, if for every affine \'etale morphism $U\to\YY$, with $U$
affine, the  fibered product $V=\XX\times_{\YY}U$ is affine.
\end{defn}

\begin{prop}
A morphism $\XX\to \YY$ is affine if there exists an affine \'etale
cover $U_i\to\YY$, such that for every $i$ the fibered
product $V_i=\XX\times_\YY U_i$ is affine.
\end{prop}
\begin{pf}
This follows from Lemmas~\ref{prodet} and~\ref{app-rep}.
\end{pf}

\begin{prop}
The diagonal $\XX\to\XX\times\XX$ of a differential graded scheme
$\XX$ is affine.
\end{prop}
\begin{pf}
If $U_i\to\XX$ is an affine \'etale cover of $\XX$, then $U_i\times
U_j$ is an affine \'etale cover of $\XX\times \XX$, and we have
2-cartesian diagrams
$$\xymatrix{
U_{ij}\dto\rto & U_i\times U_j\dto\\
\XX\rto & \XX\times \XX}$$
where $U_{ij}=U_i\times_{\XX}U_j$, which is affine, by the definition
of affine \'etale morphism.
\end{pf}

\subsection{\'Etale morphisms}

\begin{defn}
A morphism $f:\XX\to \YY$ of differential graded schemes is called
{\bf \'etale}, if for every morphism $U\to \YY$, with $U$ affine,
the fibered product $\VV=\XX\times_{\YY}U$ is a differential graded
scheme and for every affine \'etale morphism $V\to \VV$, with $V$
affine, the composition $V\to U$ is \'etale.
\end{defn}

\begin{prop}
Let $f:\XX\to \YY$ be a morphism of differential graded schemes.
Suppose given an epimorphic  family of morphisms $U_i\to\YY$, with $U_i$
affine for all $i$.  Suppose
further that for every $i$ the fibered product
$\VV_i=\XX\times_{\YY}U_i$ is a differential graded scheme and that
there exists an affine \'etale cover $V_{ij}\to\VV_i$ of $\VV_i$ such
that the composition $V_{ij}\to U_i$ is \'etale, for all $j$. 
Then $f$ is \'etale
\end{prop}
\begin{pf}
This is not difficult to prove using the techniques developed so far.
In particular, use Lemma~\ref{descjdskfl} and Note~\ref{loclet}.
\end{pf}

\begin{note}
A morphism of differential graded schemes is affine and \'etale if
and only if it is affine \'etale.  An \'etale morphism with affine
source is affine \'etale.
\end{note}

\begin{cor}
If $\XX\to\YY$ is a morphism of differential graded schemes and
$\YY'\to \YY$ an \'etale morphism of differential graded schemes, then
the fibered product $\XX'=\XX\times_\YY\YY'$ is a differential graded
scheme and $\XX'\to \XX$ is \'etale.\qed
\end{cor}

\begin{rmk}
The question of the existence of more general fibered products in the
2-category of differential graded schemes is rather subtle.  It is
treated in detail in \cite{dgsIV}.
\end{rmk}

\begin{prop}
Consider the 2-commutative diagram of differential graded schemes
$$\xymatrix@=1pc{
\XX\rrto^f\drto_h && \YY\dlto^g\\
&\ZZ &}$$
If $g$ is \'etale, then $f$ is \'etale if and only if $h$ is \'etale.
If $f$ is an \'etale epimorphism, then $g$ is \'etale if and only if
$h$ is \'etale. \qed
\end{prop}

\begin{prop}
A morphism of differential graded schemes is \'etale if and only if it
is categorically \'etale.
\end{prop}

\begin{defn}
An \'etale morphism $f:\XX\to\YY$ of differential graded schemes is
called an {\bf open immersion}, if it is a monomorphism of
differential graded sheaves.
\end{defn}

\begin{rmk}
An affine \'etale morphism is an open immersion if and only if it is
an affine open immersion.  Open immersions are stable under base change.
\end{rmk}

\begin{prop}
A morphism of differential graded schemes is an open immersion if and
only if it is a monomorphism.
\end{prop}

\begin{defn}
If $\XX$ is a differential graded scheme, an {\bf open subscheme }of
$\XX$ is a full sub-2-category $\XX'\subset \XX$, such that $\XX'$ is
itself a differential graded scheme and the inclusion morphism
$\XX'\to\XX$ is an open immersion. 
\end{defn}

\Section{The basic 1-categorical invariants}\label{sec.oneinv}

For any differential graded sheaf $\FF$, the underlying 2-category
$\FF$ is endowed with a topology in a canonical way. (A sieve is
covering in $\FF$ if its image in $\SS$ is covering.) Thus every
differential graded scheme $\XX$ has an associated 2-site, namely the
2-category $\XX$ itself, with this canonical topology.

\begin{defn}
A {\bf sheaf }over $\XX$, is a sheaf
on this associated 2-site.
\end{defn}

Write $\ol\XX$ for the 1-category associated to  $\XX$.
By the basic 1-categorical invariants of the differential graded
scheme $\XX$ we mean certain sheaves of sets on $\XX$.  Note that
every sheaf of sets on $\XX$ comes in a unique and canonical way from
a sheaf on $\ol\XX$.  Thus the terminology.

\begin{rmk}
The rule $\XX\to\ol\XX$ defines a 2-functor from the 2-category of
differential graded schemes to the 2-category of sites.
\end{rmk}

\subsection{The associated graded structure sheaf}

Let $\XX$ be a differential graded scheme. For an object $x$ of 
$\XX$, denote the image of $x$ under the structure 2-functor $\XX\to
\SS$ by $A_x$.  

\begin{defn}
The {\bf truncated structure sheaf }of $\XX$ is the sheaf of sets on
$\XX$ defined by 
$$x\longmapsto h^0(A_x)\,.$$
We denote the truncated structure sheaf by $h^0(\O_\XX)$.  This is an
abuse of notation, as we have not defined $\O_\XX$.
\end{defn}

Note that for every morphism $U\to \XX$, with $U$ affine,
$h^0(\O_\XX)(U)$ is a finitely generated $k$-algebra.

The fact that $h^0(\O_\XX)$ is a sheaf, follows directly from the
definition of the \'etale topology.

Thus $h^0(\O_\XX)$ is a sheaf of $k$-algebras on $\XX$.

\begin{defn}
The $n$-th {\bf higher structure sheaf }is the sheaf of sets on $\XX$
defined by
$$x\longmapsto h^n(A_x)\,.$$
We denote the $n$-the higher structure sheaf by $h^n(\O_\XX)$. 

The direct sum
$$h^\ast(\O_\XX)=\bigoplus_nh^n(\O_\XX)$$
is called the {\bf associated graded structure sheaf }of $\XX$.
\end{defn}

For every morphism $U\to\XX$ with $U$ affine, $h^n(\O_\XX)(U)$ is
a finitely generated $h^0(\O_\XX)(U)$-module and $h^\ast(\O_\XX)$ is a
graded $h^0(\O_\XX)$-algebra. Thus $h^n(\O_\XX)$ is a
coherent sheaf of modules and $h^\ast(\O_\XX)$ a graded sheaf of
algebras over the sheaf of $k$-algebras
$h^0(\O_\XX)$. 

Of course, $h^n(\O_\XX)=0$, for all $n>0$. 

\begin{numrmk}
Let $\phi:\XX\to\YY$ be a morphism is differential graded
schemes. Then there is a natural isomorphism ({\em sic!\/}) of sheaves of
graded
$k$-algebras
$$\phi^{-1}h\upst(\O_\YY)\longrightarrow h\upst(\O_\XX)\,.$$
For example, for $\YY=\dgspec k=\SS$, we get that $h\upst(\O_\XX)$
is the pullback of $h^\ast(\O_{\dgspec k})$ via the structure functor
$\XX\to\SS$ (which is also clear from the definition). We will
abbreviate 
$h^\ast(\O_{\dgspec k})$ by $h^\ast(\O)$. 
\end{numrmk}

\subsection{Higher tangent sheaves}

We need some preliminaries concerning the naturality properties of
$\Deru(B,A)$, for resolving algebras $B$, $A$.
(For the notation $\Deru(B,A)$, see
Section~I\ref{difcot}.)

Let $f:B\to A$ and $g:B\to A$ be morphisms of resolving algebras.  Let
$\theta:f\Rightarrow g$ be a homotopy, i.e., a morphism $\theta:B\to
A\otimes\Omega_1$, such that $\del_0\theta=f$ and $\del_1\theta=g$. 
$$\xymatrix{
B\rtwocell^f_g{_\theta} & A}$$
Recall (Definition~I.\ref{can-iso}), that $\theta$ induces a canonical
isomorphism of $h^0(B)$-modules 
$$h_\ell \Deru(B,{_fA})\stackrel{\theta\lst}{\longrightarrow}
h_\ell\Deru(B,{_gA})\,.$$ 
Recall also, that  $\theta\lst$ depends only on the homotopy class of
$\theta$, and is thus well-defined for a 
2-isomorphism $\theta:f\Rightarrow g$ in $\RR$.  Moreover, it is
functorial for vertical composition of 2-morphisms: 
$\theta\lst\eta\lst=(\theta\eta)\lst$. 

We will require two further naturality properties of this induced
canonical 
isomorphism.  First, some more notation:

Suppose given a 2-commutative diagram
$$\xymatrix{
B\drto\rto\drtwocell\omit{<-2>^\eta} & A\dto^f \\
& A'}$$
We denote the composition
$$h_\ell\Deru(B,A)\stackrel{f\lst}{\longrightarrow} h_\ell
\Deru(B,{_fA'})\stackrel{\eta\lst}{\longrightarrow} h_\ell
\Deru(B,A')$$ 
by $\eta(f\lst)$. 

\begin{prop}\label{propegg}
Consider a diagram
\vskip-1cm
\begin{equation}\label{egg}
\xymatrix{
B\rto\rruppertwocell<12>^{}{^\eta}
\rrlowertwocell<-12>_{}{_\zeta}&
A\rtwocell^f_{g}{\omit}& A'}
\end{equation}
\vskip-1cm
\noindent where the two morphisms from $B$ to $A'$ are equal. We get
two 
homomorphisms of $h^0(B)$-modules
$$\xymatrix{
{h_\ell\Deru(B,A)}\rrtwocell^{\eta(f\lst)}_{\zeta(g\lst)}{\omit}& &
{h_\ell\Deru(B,A')}}\,.$$
These are equal, if there exists a 2-isomorphism
$$\xymatrix{A\rtwocell^f_g{_\theta}&A'}$$
making Diagram (\ref{egg}) commute, i.e., such that
$\eta=\zeta\theta$.
\end{prop}
\begin{pf}
The proof is similar to the construction of the induced canonical
isomorphism. 
\end{pf}

\begin{prop}\label{tricky}
Consider a diagram
\begin{equation*}
\xymatrix{
B\rtwocell^f_{g}{_\theta} & B'\rto^p& A\,.}
\end{equation*}
Then the induced diagram
$$\xymatrix@R=.5pc{
&& {h_\ell \Deru(B,{_{f}A})}\ddto^{(p\theta)\lst}\\
{h_\ell\Deru(B',A)}\urrto^{f\upst}\drrto_{g\upst}&&\\
&& {h_\ell \Deru(B,{_{g}A})}}$$
commutes.
\end{prop}
\begin{pf}
This proposition is a little more tricky.  The problem is that
$B'\otimes\Omega_1$ is not a resolving algebra. Thus it is not clear
if $B'\otimes\Omega_1$ is cofibrant, and hence if
$\Deru(B'\otimes\Omega_1, \argument)$ is sufficiently well-behaved.  Thus,
instead of working with $\Deru(B'\otimes\Omega_1, \argument)$, we use
$\Homu_{B'\otimes\Omega_1}(L_{B'\otimes \Omega_1},\argument)$.  This
requires a theory of the cotangent complex for non-resolving algebras,
as developed in \cite{hinich}. We omit the details of the rather
lengthy diagram chases that conclude the proof.
\end{pf}

We are now ready to define the higher tangent sheaves.
Let $B$ be a perfect resolving algebra. 

\begin{defn}
The $\ell$-th {\bf higher tangent sheaf }of $\dgspec B$ is the sheaf
of sets 
on $\dgspec B$ defined by 
$$(B\stackrel{x}{\to} A)\longmapsto
h_\ell\Deru(B,A)=h_\ell(\Theta_B\otimes_BA)\,.$$ 
We denote the $\ell$-th higher tangent sheaf of $\dgspec B$ by
$h_\ell(\Theta_{\dgspec B})$. 
\end{defn}

This definition gives rise to a presheaf, because of
Proposition~\ref{propegg}. 
The fact that $h_\ell(\Theta_{\dgspec B})$ is a sheaf follows from
Corollary~I.\ref{prep.desc}.

\begin{rmk}
A morphism of perfect resolving algebras $f:B\to B'$, which gives rise
to the morphism of differential graded schemes $\phi:\dgspec
B'\to\dgspec B$, defines a canonical sheaf map
\begin{equation}\label{can-sheaf}
h_\ell(\Theta_{\dgspec B'})\longrightarrow \phi^{-1}
h_\ell(\Theta_{\dgspec B})\,.
\end{equation}
Let $\theta:f\Rightarrow g$ be a homotopy between the morphisms
$f,g:B\to B'$.  Letting $\phi$ be the morphism of differential graded
schemes induced by $f$ and $\psi$ the morphism of differential graded
schemes induced by $g$, we get an induced 2-isomorphism
$\eta:\psi\Rightarrow \phi$. The 2-isomorphism $\eta$ gives rise to a natural
equivalence of functors
$\eta^{-1}:\phi^{-1}\to \psi^{-1}$ and hence to a sheaf
isomorphism
$$\eta^{-1}:\phi^{-1}h_\ell(\Theta_{\dgspec B})\longrightarrow
\psi^{-1}h_\ell(\Theta_{\dgspec B})\,.$$
The induced triangle of sheaves
$$\xymatrix@R=.5pc{
&& {\phi^{-1}h_\ell(\Theta_{\dgspec B})}\ddto^{\eta^{-1}}\\
{h_\ell(\Theta_{\dgspec B'})}\urrto\drrto &&\\
&& {\psi^{-1}h_\ell(\Theta_{\dgspec B})}}$$
commutes.  This follows from Proposition~\ref{tricky}.

Because of this, we may define $h_\ell(\Theta_U)$ for any affine
differential graded scheme $U$, because the sheaf on $U$ pulled back
via any isomorphism $U\to\dgspec B$, is independent of the choice of
$B$ and $U\to\dgspec B$, at least up to canonical isomorphism. 
\end{rmk}

Let $\XX$ be a differential graded scheme.  Let $V\to\XX$ and
$U\to\XX$ be \'etale morphisms, with $V$ and $U$ affine.  Assume given
a 2-commutative diagram of differential graded schemes
\begin{equation}\label{afetx}
\vcenter{\xymatrix{
V\rto^f\drto\drtwocell\omit{<-2>_{}} & U\dto \\ & \XX}}
\end{equation}
we get an induced sheaf map
$$h_\ell(\Theta_V)\longrightarrow f^{-1}h_\ell(\Theta_U)\,,$$
which is an isomorphism, by Proposition~I.\ref{mea}, and because $V\to
U$ is necessarily \'etale.  

Thus, as we let $U\to\XX$ vary over all \'etale morphisms with affine
$U$, we get gluing data for a sheaf of sets $h_\ell(\Theta_\XX)$ on
$\XX$.

\begin{defn}
The sheaf $h_\ell(\Theta_\XX)$  is called the $\ell$-th {\bf higher
tangent sheaf }of $\XX$.
The {\bf associated graded tangent sheaf }of $\XX$ is the direct sum
$$h_\ast(\Theta_\XX)=\bigoplus_\ell h_\ell(\Theta_\XX)\,.$$
\end{defn}

The higher tangent sheaf $h_\ell(\Theta_\XX)$ comes with isomorphisms
$$h_\ell(\Theta_U)\longrightarrow h_\ell(\Theta_\XX)\resto U\,,$$
for every \'etale $U\to \XX$ with affine $U$. Any diagram
(\ref{afetx}) induces a commutative diagram
$$\xymatrix{
h_\ell(\Theta_V)\dto\rto & f^{-1} h_\ell(\Theta_U)\dto\\
{h_\ell(\Theta_\XX)\resto V} \rto &
f^{-1}\big(h_\ell(\Theta_\XX)\resto U\big) }$$
of sheaves on $V$.

Every higher tangent sheaf $h_\ell(\Theta_{\XX})$ is a coherent
$h^0(\O_\XX)$-module.  The associated graded tangent sheaf
$h_\ast(\Theta_{\XX})$ is a sheaf of graded
$h^\ast(\O_\XX)$-modules.

Note that $h_\ell(\Theta_\XX)=0$, for
$\ell<-N$, if $\XX$ is of amplitude $N$.

\begin{rmk}
If $\dgspec B\to \XX$ is \'etale and $B\to A$ is an
arbitrary morphism of perfect resolving algebras, then we have, by
construction of $h_\ell(\Theta_\XX)$, a canonical isomorphism of
$h^0(B)$-modules
$$h_\ell\Deru(B,A)\longiso h_\ell(\Theta_\XX) (\dgspec A)\,.$$
For example, if $\phi:\dgspec K\to \XX$ is  a $K$-valued point of $\XX$,
then 
$$h_\ell(\Theta_\XX)(\phi)=h_\ell\Deru(B,K)\,,$$
for any affine \'etale neighbourhood $\dgspec K\to\dgspec B\to\XX$ of
$\phi$. 
\end{rmk}

\begin{defn}
The {\bf higher tangent spaces }of $\XX$ at the $K$-valued point
$\phi$ of $\XX$ are the $K$-vector spaces
$h_\ell(\Theta_\XX)(\phi)$. 
\end{defn}

\begin{prop}
If $f:\XX\to\YY$ is a morphism of differential graded schemes, then we
get an induced canonical homomorphism of sheaves of
$h^0(\O_\XX)$-modules
\begin{equation}\label{cantheta}
h_\ell(\Theta_\XX)\longrightarrow f^{-1} h_\ell(\Theta_\YY)\,.
\end{equation}
The morphism  $f$ is \'etale, if and only if 
(\ref{cantheta}) is an isomorphism, for all $\ell$ (or for one fixed
$\ell\geq0$). 

If $\phi:\dgspec K\to\XX$ is a $K$-valued point of $\XX$, then $f$
induces a canonical homomorphism of $K$-vector spaces
\begin{equation}\label{cannottheta}
h_\ell(\Theta_\XX)(\phi)\longrightarrow
h_\ell(\Theta_\YY)\big(f(\phi)\big)\,,
\end{equation}
for all $\ell$. 
The morphism $f$ is \'etale, if and only if (\ref{cannottheta})  is an
isomorphism, for all $\ell$.\qed
\end{prop}

\subsubsection{The relative case}

Note that if $C'\to C\to B\to A$ is a composition of morphisms of
resolving algebras, with $C$, $C'$ and $B$ perfect and $C\to B$ as
well as $C'\to B$ resolving, then the canonical homomorphism of
$h^0(A)$-modules
$$h_\ell\Deru_C(B,A)\longrightarrow h_\ell\Deru_{C'}(B,A)$$
is an isomorphism, if $C'\to C$ is \'etale. 

Now assume given a morphism of differential graded schemes
$\phi:\XX\to\YY$. 

\begin{defn}
The $\ell$-th {\bf relative higher tangent sheaf }of $\XX\to\YY$,
notation $h_\ell(\Theta_{\XX/\YY})$, is defined in such a way that for
every resolving morphism of perfect resolving algebras $C\to B$ and any
2-commutative diagram of differential graded schemes
\begin{equation}\label{reltan}
\vcenter{
\xymatrix{
{\dgspec B}\dto\rto \drtwocell\omit{^} & {\dgspec C}\dto\\
{\XX}\rto & {\YY}}}
\end{equation}
where $\dgspec B\to\XX$ and $\dgspec C\to\YY$ are \'etale, we have a
canonical isomorphism
$$h_\ell(\Theta_{B/C})\longrightarrow
h_\ell(\Theta_{\XX/\YY})\resto\dgspec B\,,$$
where $h_\ell(\Theta_{B/C})$ is the sheaf on $\dgspec B$ defined by 
$$(B\to A)\longmapsto h_\ell \Deru_C(B,A)\,.$$
\end{defn}

By the remark preceding the definition, $h_\ell(\Theta_{B/C})$ does
not depend on the choice of $C$. Note also that the \'etale morphisms
$\dgspec B\to \XX$ admitting a factorization (\ref{reltan}) are
cofinal in the 2-category of all \'etale $\dgspec B\to\XX$ and still
cover $\XX$.  Thus the fact that $h_\ell(\Theta_{\XX/\YY})$ exists
with the required properties is proved as in the absolute case.

\begin{numrmk}
The morphism $\XX\to \YY$ is \'etale, if and only if
$h_\ell(\XX/\YY)=0$ for all $\ell$, or for on fixed $\ell\geq0$. 
\end{numrmk}

\begin{numrmk}
Given $C\to B$ and a diagram such as (\ref{reltan}), for every
morphism of resolving algebras $B\to A$ we have a short exact sequence
of complexes of $A$-modules
$$\ses{\Deru_C(B,A)}{}{\Deru(B,A)}{}{\Deru(C,A)}\,.$$
These give rise to natural long exact sequences of
$h^0(\O_\XX)$-modules
\begin{equation}\label{grnle}
\ldots \longrightarrow
h_\ell(\Theta_{\XX/\YY}) \longrightarrow
h_\ell(\Theta_\XX)\longrightarrow
\phi^{-1}h_\ell(\Theta_\YY)\longrightarrow
h_{\ell-1}(\Theta_{\XX/\YY})\longrightarrow \ldots
\end{equation}

Let 
\begin{equation}\label{justtwo}
\vcenter{\xymatrix{
{\XX}\rto^\phi\drto_\kappa\drtwocell\omit{^<-2>} & {\YY}\dto^\psi \\
&{\ZZ}}}
\end{equation}
be a 2-commutative diagram of differential graded schemes.  By the
naturality properties of~(\ref{grnle}), we have a `long exact braid
with four strands'
\vskip1pc
\begin{equation}\label{braid}
\vcenter{
\xymatrix@C=4pc@R=3pc@!0{
& {\phi^{-1}h_\ell(\Theta_{\YY/\ZZ})}\ar@/^2pc/[rr]\drto &&
{h_{\ell-1}(\Theta_{\XX/\YY})}\ar@/^2pc/[rr]\drto && 
{h_{\ell-1}(\Theta_{\XX})} \\
{h_\ell(\Theta_{\XX/\ZZ})}\urto\drto &&
{\phi^{-1}h_{\ell}(\Theta_{\YY})}\urto\drto  &&
{h_{\ell-1}(\Theta_{\XX/\ZZ})}\urto\drto &\\
& {h_\ell(\Theta_{\XX})}\ar@/_2pc/[rr]\urto &&
{\kappa^{-1}h_{\ell}(\Theta_{\ZZ})}\ar@/_2pc/[rr]\urto && 
{\phi^{-1}h_{\ell-1}(\Theta_{\YY/\ZZ})} 
}}\end{equation}
\vskip1cm
\end{numrmk}

\subsection{Homotopy sheaves}

The homotopy sheaves of a differential graded scheme are defined
similarly to the higher tangent sheaves.  First, they are defined for
affine differential graded schemes and then they are glued with
respect to the \'etale topology.

Let $B$ be a perfect resolving algebra and $\ell>0$ and integer.

\begin{defn}
The $\ell$-th {\bf homotopy sheaf }of $\dgspec B$, notation
$\pi_\ell(\dgspec B)$, is the sheaf of
sets on $\dgspec B$ defined by
$$(B\stackrel{x}{\to} A)\longmapsto
\pi_\ell\shom(B,A)\,.$$
\end{defn}

The fact that this defines a sheaf on $\dgspec B$ follows directly
from descent theory, Theorem~\ref{Descent}(i). 

A morphism of perfect resolving algebras $f:B\to B'$ gives rise to a
canonical sheaf map
\begin{equation}\label{seccan}
\pi_\ell(\dgspec B')\longrightarrow \phi^{-1}\pi_\ell(\dgspec
B)\,,
\end{equation}
via the restriction map $\shom(B',A)\to\shom(B,A)$. For \'etale $f$,
the sheaf map (\ref{seccan}) is an isomorphism, by
Proposition~I.\ref{etahom}.  (Here we use the fact that $\ell>0$.) 

Let now $\XX$ be a differential graded scheme. 
As we let $\dgspec B\to\XX$ vary over all perfect resolving algebras $B$
and all \'etale morphisms to $\XX$, the various $\pi_\ell(\dgspec B)$
glue, via the gluing maps~(\ref{seccan}), to a sheaf of sets $\pi_\ell(\XX)$
on $\XX$.

\begin{defn}
The sheaf $\pi_\ell(\XX)$ is called the $\ell$-th {\bf homotopy sheaf
}of the differential graded scheme $\XX$.
\end{defn}

By construction, $\pi_\ell(\XX)$ is endowed with a sheaf isomorphism
$$\pi_\ell(\dgspec B)\longrightarrow \pi_\ell(\XX)\resto \dgspec
B\,,$$ for all \'etale $\dgspec B\to\XX$. This isomorphism is
compatible with the gluing isomorphisms~(\ref{seccan}).  Thus, if $A$
is an arbitrary resolving algebra in $\RR$, endowed with a morphism
$B\to A$, then
$$\pi_\ell(\XX)(A)=\pi_\ell\shom(B,A)\,.$$

The sheaves $\pi_\ell(\XX)$ are sheaves of groups; abelian,
for $\ell\geq2$. 

\begin{numrmk}\label{autpione}
Let $\AAut\XX$ denote the sheaf of sets on $\XX$ given by 
$$x\longmapsto \Aut(x)\,.$$
Here $\Aut(x)$ stands for the automorphism group of the object $x$ of
$\XX$ inside the fiber $\XX_A$, where $x$ lies over $A\in\RR$.
For $\XX=\dgspec B$, we have 
$\pi_1(\XX)(x)=\Aut(x)$, for every object $x:B\to A$ of $\dgspec
B$. Hence $\pi_1(\XX)=\Aut\XX$.

For every morphism of differential graded schemes $\phi:\YY\to\XX$, we
have an induced morphism $\AAut\YY\to\phi^{-1}\AAut\XX$ of sheaves of
sets on $\YY$.  In particular, for a morphism $\phi:\dgspec B\to \XX$,
we get a canonical morphism $\pi_1(\dgspec B)\to \phi^{-1}\AAut\XX$.
These canonical morphisms glue to give a canonical morphism
$$\pi_1(\XX)\longrightarrow\AAut\XX\,,$$
which is trivially an isomorphism.

If we choose pullbacks for the fibered category $\XX\to\SS$, we can
identify the fiber $\XX_A$ of $\XX$ over the differential graded
algebra $A$ in $\RR$, with the groupoid $\bhom(\dgspec A,\XX)$. Doing
this we have
$$\pi_1(\XX)(x)=\pi_1\bhom(\dgspec A,\XX)\,,$$
for any $x:\dgspec A\to \XX$.
\end{numrmk}

\begin{defn}
Let $\pi_0(\XX)$ denote the presheaf of pointed sets on $\XX$ defined by 
$$\pi_0(\XX)(x)=\pi_0(\XX_{A_x})\,,$$
where $A_x$ is the image of $x$ in $\SS$. Note that by
Theorem~\ref{Descent}~(i), for affine $\XX$, the presheaf $\pi_0(\XX)$
is a sheaf.
\end{defn}

\begin{example}
Denote $\dgspec k[x]$, where $\deg x=0$, by $\aaa^1$. Then we have
$\pi_\ell(\aaa^1)= h^{-\ell}(\O_{\aaa^1})$. Thus, we
have that $\pi_\ell(\aaa^1)$ is the pullback of $h^{-\ell}(\O)$ via
the structure functor $\aaa^1\to\SS$, for all $\ell\geq0$. 

Let us call a morphism of differential graded schemes
$\phi:\XX\to\aaa^1$ a {\bf regular function }on $\XX$. Then for any
regular function $\phi$ on the differential graded scheme $\XX$ we
have
$$h^{-\ell}(\O_\XX)=\phi^{-1}\pi_\ell(\aaa^1)\,,$$
for all $\ell\geq0$. 

We have a canonical map
$$\pi_0\Homu(\XX,\aaa^1)\longrightarrow
\Gamma\big(\XX,h^0(\O_\XX)\big)\,,$$
which is bijective, if $\XX$ is affine.
\end{example}

\subsubsection{The relative case}

Just like the higher tangent sheaves, the homotopy sheaves also admit
relative versions. Let $\phi:\XX\to \YY$ be a morphism of
differential graded schemes and $\ell>0$ an integer. 

\begin{defn}
The $\ell$-th {\bf relative homotopy sheaf }of $\XX$ over $\YY$,
notation $\pi_\ell(\XX/\YY)$,  is
defined in such a way that for every resolving morphism of perfect
resolving algebras $C\to B$ and any 2-commutative diagram of
differential graded schemes (\ref{reltan}), where $\dgspec B\to \XX$
and $\dgspec C\to \YY$ are \'etale, we have a conical isomorphism 
$$\pi_\ell(B/C)\longrightarrow\pi_\ell(\XX/\YY)\resto\dgspec B\,,$$
where $\pi_\ell(B/C)$ is the sheaf on $\dgspec B$ defined by 
$$(B\to A)\longmapsto \pi_\ell\shom_C(B,A)\,.$$

Moreover, define the presheaf of pointed sets on $\XX$ 
$$\pi_0(\XX/\YY)$$
by defining $\pi_0(\XX/\YY)(x)$ to be $\pi_0$ of the fiber through $x$
of the morphism of groupoids $\XX_{A}\to\YY_{A}$, where $A$ is
the object of $\RR$ over which the object $x$ of $\XX$ lies.
\end{defn}

\begin{prop}
Let $\phi:\XX\to\YY$ be a morphism of differential graded schemes and
$r>0$ and integer. The
following are equivalent:

(i) $\pi_\ell(\XX/\YY)=0$, for all $\ell>0$,

(ii) $\pi_r(\XX/\YY)=0$,

(iii) $\pi_\ell(\XX)\to\phi^{-1}\pi_\ell(\YY)$
is an isomorphism of sheaves of groups, for all $\ell>0$,

(iv) $\pi_r(\XX)\to\phi^{-1}\pi_r(\YY)$
is an isomorphism of sheaves of groups.
\end{prop}
\begin{pf}
Follows from Proposition~I.\ref{etahom}.
\end{pf}

\begin{prop}[Long exact homotopy sequence]
There is a natural long sequence of presheaves on
$\XX$ 
\begin{multline}\label{lesp}
\ldots\longrightarrow \pi_\ell(\XX/\YY) \longrightarrow \pi_\ell(\XX)
\longrightarrow \phi^{-1}\pi_\ell(\YY)\stackrel{\del}{\longrightarrow}
\pi_{\ell-1}(\XX/\YY)\longrightarrow\ldots\\
\ldots\longrightarrow \phi^{-1}\pi_1(\YY)\stackrel{\del}{\longrightarrow}
\pi_0(\XX/\YY) \longrightarrow \pi_0(\XX)\longrightarrow \pi_0(\YY)\,.
\end{multline}
This sequence gives rise to a long exact sequence of groups and
pointed sets, 
when evaluated at an 
object $x:\dgspec A\to \XX$ of $\XX$, which admits a factorization
$$
\xymatrix{
{\dgspec A}\dto_x\rto \drtwocell\omit{^} & {\dgspec C}\dto\\
{\XX}\rto_\phi & {\YY}}
$$
with \'etale $\dgspec C\to \YY$.  In particular, the part of
(\ref{lesp})  ending with $\phi^{-1}\pi_1(\YY)$ is an exact sequence of
sheaves of groups on $\XX$.
\end{prop}
\begin{pf}
Let $C\to B$ be a resolving morphism of perfect resolving algebras
together with a 2-commutative diagram 
\begin{equation}\label{phipsi}
\vcenter{
\xymatrix{
{\dgspec B}\dto\rto^\psi \drtwocell\omit{^} & {\dgspec C}\dto\\
{\XX}\rto_\phi & {\YY}}}
\end{equation}
with $\dgspec B\to
\XX$ and $\dgspec C\to \YY$ \'etale. Then for any morphism of
resolving algebras $B\to A$, we get a fibration of spaces
$$\shom(B,A)\longrightarrow\shom(C,A)\,$$
with fiber $\shom_C(B,A)$. The associated long exact homotopy
sequence, which is natural in $A$, gives rise to a long exact sequence
of presheaves on $\dgspec B$
$$
\ldots\longrightarrow \pi_\ell(B/C) \longrightarrow \pi_\ell(B)
\longrightarrow \psi^{-1}\pi_\ell(C)\longrightarrow\ldots\longrightarrow
\psi^{-1}\pi_1(C)\,.
$$
One checks that the various maps in this sequence glue, to give a
sequence of sheaves on $\XX$
\begin{equation}\label{almost}
\ldots\longrightarrow \pi_\ell(\XX/\YY) \longrightarrow \pi_\ell(\XX)
\longrightarrow \phi^{-1}\pi_\ell(\YY)\longrightarrow\ldots\longrightarrow
\phi^{-1}\pi_1(\YY)\,.
\end{equation}
By construction, this sequence is exact on the level of groups of
sections over any object $x$ of $\XX$, lying over $A$ in $\RR$, such
that $x:\dgspec A\to \XX$ factors through a diagram~(\ref{phipsi}). 

By construction, we have, for any object $x$ of $\XX$, an exact
sequence of groups and pointed sets
$$\Aut(x)\longrightarrow \Aut\big(\phi(x)\big)\longrightarrow
\pi_0(\XX/\YY)(x)\longrightarrow \pi_0(\XX)(x)\longrightarrow 
\pi_0(\YY)\big(\phi(x)\big)\,.$$
(This is a general fact about morphisms of groupoids.)
Thus we get an exact sequence of presheaves of groups and pointed sets
on $\XX$
$$\AAut(\XX)\longrightarrow \phi^{-1}\AAut(\YY)\longrightarrow
\pi_0(\XX/\YY)\longrightarrow \pi_0(\XX)\longrightarrow
\phi^{-1}\pi_0(\YY)\,. $$
By Remark~\ref{autpione}, we have natural identifications
$\pi_1(\XX)=\AAut(\XX)$ and $\pi_1(\YY)=\AAut(\YY)$, and so we can
extend the sequence (\ref{almost}) three steps further to the right,
as required.
\end{pf}

\begin{numrmk}
Given a 2-commutative diagram of differential graded
schemes~(\ref{justtwo}), then the various long exact homotopy
sequences (\ref{lesp}) are natural enough to give rise to a
commutative long exact braid with four strands, similar
to~(\ref{braid}), except that it has a right end.
\end{numrmk}

\subsection{Differentials}

Let $\phi:\XX\to \YY$ be a morphism of differential graded
schemes. 

\begin{prop}
There exists an object $\ol\Omega_{\XX/\YY}$ in the derived category
of $h^0(\O_\XX)$, together with natural isomorphisms 
$$\ol\Omega_{\XX/\YY}\resto\dgspec B\longiso \ol\Omega_{B/C}\,,$$
for every resolving morphism of perfect resolving algebras $C\to B$
and every 2-commutative diagram (\ref{phipsi}) with \'etale $\dgspec
B\to \XX$ and $\dgspec C\to \YY$. Here $\ol\Omega_{B/C}$ is the
complex of sheaves of (finitely generated, free) $h^0(\O_{\dgspec
B})$-modules defined by 
$$(B\stackrel{x}{\to} A)\longmapsto
\ol\Omega_{B/C}(x)=\Omega_{B/C}\otimes_Bh^0(A)\,.$$ 
The complex $\ol\Omega_{\XX/\YY}$ is perfect.
\end{prop}
\begin{pf}
Note that
$\Omega_{B/C}\otimes_Bh^0(A)=
\Omega_{B/C}\otimes_Bh^0(B)\otimes_{h^0(B)}h^0(A)$. 
Hence $\ol\Omega_{B/C}$ is a presheaf of complexes of $h^0(\O_{\dgspec
B})$-modules. Because every one of the complexes of sections of
$\ol\Omega_{B/C}$ is finitely generated and free, $\ol\Omega_{B/C}$
is, in fact, a complex of sheaves of $h^0(\O_{\dgspec B})$-modules.
Thus we have constructed an object $\ol\Omega_{B/C}$ in the derived
category of $h^0(\O_{\dgspec B})$-modules which are bounded above and
have coherent cohomology.

This construction gives rise to a functor from the (1-category
associated to the) category of all
$C\to B$ and Diagrams~(\ref{phipsi}), to the derived category
$h^0(\O_\XX)$-modules, and by cohomological descent, we get the
required object 
$\ol\Omega_{\XX/\YY}$. 
\end{pf}

\begin{defn}
The complex $\ol\Omega_{\XX/\YY}$ is called the {\bf cotangent complex
} or the {\bf complex of differentials }of $\XX$ over $\YY$. Its dual
is denoted by $\ol\Theta_{\XX/\YY}$ and is called the {\bf tangent
complex }of $\XX$ over $\YY$.
\end{defn}

\begin{prop}\label{inddistri}
Given a 2-commutative diagram (\ref{justtwo}) of differential graded
schemes, we get induced distinguished triangles
$$\phi^{-1}\ol\Omega_{\YY/\ZZ}\longrightarrow
\ol\Omega_{\XX/\ZZ}\longrightarrow \ol\Omega_{\XX/\YY}\longrightarrow
\phi^{-1}\ol\Omega_{\YY/\ZZ}[1]$$
and
$$\ol\Theta_{\XX/\YY}\longrightarrow
\ol\Theta_{\XX/\ZZ}\longrightarrow
\phi^{-1}\ol\Theta_{\YY/\ZZ}\longrightarrow 
\ol\Theta_{\XX/\YY}[1]$$
in the derived category of $h^0(\O_\XX)$. These distinguished triangles
are natural in the sense that they  give rise to  commutative
`octahedra'. Let us only display the octahedron for $\ol\Theta$:
\begin{equation*}
\vcenter{
\xymatrix@C=4pc@R=3pc@!0{
&& {\phi^{-1}\ol\Theta_{\YY/\ZZ}}\ar@/^2pc/[rr]\drto &&
{\ol\Theta_{\XX/\YY}[1]}\drto &&  \\
&{\ol\Theta_{\XX/\ZZ}}\urto\drto &&
{\phi^{-1}\ol\Theta_{\YY}}\urto\drto  &&
{\ol\Theta_{\XX/\ZZ}[1]}\drto &\\
{\ol\Theta_{\XX/\YY}}\urto\ar@/_2pc/[rr]&&
{\ol\Theta_{\XX}}\ar@/_2pc/[rr]\urto && 
{\kappa^{-1}\ol\Theta_{\ZZ}}\ar@/_2pc/[rr]\urto && 
{\phi^{-1}\ol\Theta_{\YY/\ZZ}[1]} 
}}\end{equation*}
\vskip.5cm
\end{prop}

\begin{defn}
If $\ol\Omega_{\XX/\YY}$ has perfect amplitude contained in $[-N,0]$,
we say that $\phi:\XX\to\YY$ has {\bf amplitude }$N$. If $\phi:\XX\to\YY$
has amplitude $N$, we write $N=\amp(\XX/\YY)$. 
\end{defn}

Note that this definition of amplitude agrees with the earlier one for
the absolute case, Definition~\ref{def-amp}.

\begin{cor}
We have
\begin{align*}
\amp(\XX)&=\max\big(\amp(\YY),\amp(\XX/\YY)\big)\,,\\
\amp(\XX/\YY)&=\max\big(\amp(\YY)+1,\amp(\XX)\big)\,,\\
\amp(\YY)&=\max\big(\amp(\XX/\YY)-1,\amp(\XX)\big)\,.
\end{align*}
\end{cor}
\begin{pf}
This follows from Proposition~\ref{inddistri}, see also
Remark~I.\ref{formnp}. Note that we also have relative versions of 
these statements, with respect to a composition $\XX\to\YY\to\ZZ$.
\end{pf}

\begin{prop}[Spectral sequence]
There is a natural convergent third quadrant spectral sequence of
coherent $h^0(\O_\XX)$-modules
$$E_2^{p,q}=h^q(\O_\XX)\otimes_{h^0(\O_\XX)}h^p(\ol\Theta_{\XX/\YY})
\Longrightarrow h^{p+q}(\Theta_{\XX/\YY})\,.$$
If $\XX\to\YY$ has amplitude $N$, then all terms of this spectral
sequence with $p+q>N$ vanish.
\end{prop}
\begin{pf}
Glue the spectral sequences from Proposition~I.\ref{ss.th} together.
\end{pf}

\subsubsection{Locally free morphisms}

\begin{defn}
Let $\phi:\XX\to\YY$ be a morphism of differential graded
schemes. Suppose there  exist affine \'etale covers $\dgspec B_i\to \XX$
and $\dgspec C_i\to\YY$, resolving morphisms $C_i\to B_i$ and
2-commutative diagrams
$$\xymatrix{
{\dgspec B_i}\dto\rto\drtwocell\omit{^} & {\dgspec C_i}\dto \\
{\XX}\rto & {\YY}}$$
where for every $i$, there exists a basis $(x_\nu)_{\nu\in I_i}$ for
$B_i$ over $C_i$, such that $dx_\nu\in C_i$, for all $\nu\in I_i$. 
In this case we call $\phi$ {\bf locally free}.
\end{defn}

Note that for a locally free morphism of differential graded schemes
$\XX\to \YY$, the cotangent complex $\ol\Omega_{\XX/\YY}$ has
locally free cohomology sheaves over $h^0(\O_\XX)$, and hence is
locally isomorphic (in the derived category of $h^0(\XX)$) to a finite
complex of finitely generated free modules with zero differential.

\begin{example}
Every \'etale morphism if locally free.
\end{example}

\begin{prop}
For every morphism $\XX\to \YY$ of  differential graded schemes
there exists an \'etale cover $\XX_i\to \XX$, such that each
composition $\XX_i\to\YY$ factors into finitely
many locally free morphisms.
\end{prop}
\begin{pf}
Use Proposition~\ref{dglfin}.
\end{pf}

\begin{prop}
There are natural  isomorphisms of sheaves of sets on $\XX$
$$\Xi_\ell:h_\ell(\Theta_{\XX/\YY})\longrightarrow\pi_\ell(\XX/\YY)\,,$$ 
for all $\ell>0$. For $\ell\geq2$, these are isomorphisms of sheaves of
abelian groups. 
If $\phi$ is locally free, then  $\Xi_1$ is also an isomorphism of
sheaves of groups.

To make the naturality properties of $\Xi_\ell$ more precise, let 
\begin{equation*}
\vcenter{\xymatrix{
{\XX}\rto^\phi\drto_\kappa\drtwocell\omit{^<-2>} & {\YY}\dto^\psi \\
&{\ZZ}}}
\end{equation*}
be a 2-commutative diagram of differential graded schemes.  Then we
have induced commutative diagrams
$$\xymatrix{
{h_\ell(\Theta_{\XX/\YY})}\rto\dto_{\Xi_\ell} &
{h_\ell(\Theta_{\XX/\ZZ})}\rto\dto^{\Xi_\ell} &
{\phi^{-1}h_\ell(\Theta_{\YY/\ZZ})}\dto^{\Xi_\ell} \\
{\pi_\ell({\XX/\YY})}\rto &
{\pi_\ell({\XX/\ZZ})}\rto &
{\phi^{-1}\pi_\ell({\YY/\ZZ})}}$$
If, moreover, $\phi$ is locally free, then we also have a commutative
diagram
$$\xymatrix{
{\phi^{-1}h_\ell(\Theta_{\YY/\ZZ})}\dto_{\Xi_\ell}\rto^\delta &
{h_{\ell-1}(\Theta_{\XX/\YY})}\dto^{\Xi_{\ell-1}} \\
{\phi^{-1}\pi_\ell({\YY/\ZZ})}\rto^\del &
{\pi_{\ell-1}({\XX/\YY})}}$$
in other words, we get a homomorphism from the long exact sequence of
higher relative tangent sheaves to the long exact sequence of relative
homotopy 
sheaves, both truncated at the transition from $\ell=1$ to $\ell=0$. 
\end{prop}
\begin{pf}
This follows from the results of Section~I.\ref{sec.lin}.
\end{pf}

\subsection{The associated algebraic space}

\begin{prop}\label{assalg}
There exists a natural 2-functor
$$h^0:(\text{\rm dg-schemes})\longrightarrow(\text{\rm algebraic 
spaces})$$
from the 2-category of differential graded schemes to the 1-category
of algebraic spaces over $k$, which satisfies the following
properties:

(i) For every perfect resolving algebra $B$ we have
$$h^0(\dgspec B)=\spec h^0(B)\,,$$

(ii) \'etale morphisms get mapped to \'etale morphisms,

(iii) any 2-cartesian diagram in which the two vertical morphisms are
\'etale 
is mapped to a 1-cartesian diagram (with two vertical \'etale
morphisms),

(iv) any affine \'etale cover gets mapped to an affine \'etale cover.

Moreover, for every differential graded scheme $\XX$, the
algebraic space $h^0(\XX)$ is locally of finite type 
over $k$ and has affine diagonal.
\end{prop}

\begin{defn}
We call $h^0(\XX)$ the algebraic space {\bf associated }to the
differential graded scheme $\XX$, or the {\bf truncation }of $\XX$.
\end{defn}

Let $\XX$ be a differential graded scheme and $\ol\XX$ its associated
1-site. Let $X=h^0(\XX)$ be the associated algebraic space, which we
consider as a fibered category over the category of finite type
affine $k$-schemes. Thus $X$ has an induced \'etale topology. Then we
can use Proposition~\ref{assalg} to 
construct a natural functor 
$$h^0:\ol\XX\longrightarrow X\,,$$
which fits into a commutative diagram of categories
$$\xymatrix{
\ol\XX\dto\rto^{h^0} & X\dto\\
\SS\rto^-{h^0}& (\text{\rm finite type affine $k$-schemes})}
$$
The functor $h^0:\ol\XX\to X$ is a continuous functor of categories
endowed with 
Grothendieck topologies. Thus sheaves
on $X$ pull back via $h^0$ to sheaves on $\ol\XX$.  Let us denote this
pull back functor by $\iota\lst$.  It has a left adjoint $\iota^{-1}$,
which extends $h^0$ from representable sheaves to all sheaves:
$$\xymatrix{
(\text{\rm sheaves on $\ol\XX$})\ar@<.5ex>[r]^{\iota^{-1}} &
(\text{\rm sheaves on $X$})\ar@<+.5ex>[l]^{\iota\lst}}$$
Thus we have a morphism of topoi
$$\iota:(\text{\rm sheaves on $X$})\longrightarrow(\text{\rm sheaves
on $\ol\XX$}) $$
and a morphism of sites
$$\iota:X\longrightarrow\ol\XX\,.$$
This morphism of sites should be thought of as a globalization and a
dualization of the natural morphism of differential graded algebras
$A\to h^0(A)$. But note that we cannot think of the algebraic space
$X$ as a differential graded scheme, unless it has  perfect cotangent
complex. (See \cite{dgsIII}.)

Note that $\iota^{-1} h^0(\O_\XX)=\O_X$.  Because $\iota^{-1}$ is
exact, we can pull back  the higher structure sheaves, higher tangent sheaves,
and the tangent and cotangent complex from $\ol\XX$
to $X$, simply by applying $\iota^{-1}$. All the exactness properties
are preserved under this operation.  We will use the notation
$\argument\otimes\O_X$, to denote the functor $\iota^{-1}$. For
example,
$$\ol\Omega_{\XX/\YY}\otimes \O_X=\iota^{-1}\ol\Omega_{\XX/\YY}\,.$$

\begin{prop}
A morphism $\XX\to\YY$ of differential graded schemes if \'etale, if
and only if $\ol\Omega_{\XX/\YY}$ is acyclic and if and only if
$\ol\Omega_{\XX/\YY}\otimes \O_X$ is acyclic. \qed
\end{prop}

\begin{prop}
There is a one-to-one correspondence between the open subschemes of a
differential graded scheme $\XX$ and the open subspaces of the
associated algebraic space $X$.\qed
\end{prop}

Define a {\em closed point }of the differential graded scheme $\XX$ to
be an equivalence class of 
morphisms $\dgspec K\to\XX$, where $K$ is a finite extension field of
$k$. The equivalence relation is generated by considering $\dgspec
K\to \XX$ and 
$\dgspec L\to \XX$ equivalent, if there exists a 2-commutative diagram
$$\xymatrix{
{\dgspec K}\rto\drto\drtwocell\omit{^<-2>} & {\dgspec L}\dto\\
& \XX}$$
Let $|\XX|$ denote the set of all closed points of $\XX$.  It is a
topological space by calling a subset open if it is the set of all
closed points of an open subscheme of $\XX$.  We call $|\XX|$ the {\em
Zariski topological space }associated to $\XX$. 

If $|X|$ is the set of closed points of the algebraic space $X$
associated to $\XX$, endowed with the Zariski topology, then there is
a homeomorphism $|X|\to|\XX|$. 

\begin{defn}
For a morphism of differential graded schemes $f:\XX\to\YY$, we call
$\rk\ol\Omega_{\XX/\YY}$ the {\bf relative dimension }of $\XX$ over
$\YY$. 
\end{defn}

The relative dimension is a locally constant, integer-valued function
on $|\XX|$. 

\begin{prop}
A morphism of differential graded schemes is an
isomorphism, if and only if it is \'etale and induces an isomorphism
on truncations.
\end{prop}
\begin{pf}
This follows from Corollary~I.\ref{qiscondition}.
\end{pf}

\subsubsection{Obstruction theory}

\begin{prop}
Let $\XX\to \YY$ be a morphism of differential graded schemes and
$X\to Y$ its truncation.  Then  we have a canonical morphism in the
derived category of $\O_X$
$$\alpha:\ol\Omega_{\XX/\YY}\otimes\O_X\longrightarrow L_{X/Y}\,,$$
where $L_{X/Y}$ is the relative cotangent complex of the morphism of
algebraic spaces $X\to Y$.  The morphism $\alpha$ is a relative obstruction 
theory for $X$ over $Y$ in the sense of \cite{BF}. \qed
\end{prop}

Let $f:\XX\to\YY$ be a morphism of differential graded schemes of
amplitude 1 and $X\to
Y$ its truncation.  Let $d$ be the relative dimension of $\XX$ over
$\YY$. For every pullback diagram of algebraic
spaces
$$\xymatrix{
U\rto\dto & V\dto\\
X\rto & Y}$$
we get an induced relative obstruction theory
$\ol\Omega_{\XX/\YY}\otimes\O_U\to L_{U/V}$ for $U$ over $V$, which is
perfect, in the terminology of~\cite{BF}. 

If $V$ is a variety, then this perfect obstruction theory defines a
virtual fundamental class $f\upsh[V]\in A\lst(U)$. As we let $V$ vary,
we get a bivariant class $f\upsh\in A_d(X\to Y)$, or, in other words,
an {\em orientation }of $X$ over $Y$. 

\begin{prop}
Let $f:\XX\to \YY$ and $g:\YY\to\ZZ$ be morphisms of differential
graded schemes, both of amplitude 1. Let $h:\XX\to\ZZ$ be isomorphic
to the composition $g\comp f$. Then $h$ is also of amplitude 1 and
$h\upsh=g\upsh \cdot f\upsh$. \qed
\end{prop}

In particular, any differential graded scheme $\XX$ of amplitude 1 has
a virtual fundamental class $[\XX]\in A_d(\XX)$, where $d=\dim\XX$.
If $f:\XX\to\YY$ is a morphism of amplitude 1 between differential
graded schemes of amplitude 1, then we have $f\upsh[\YY]=[\XX]$.

\end{document}